\documentclass[journal]{IEEEtran}

 \usepackage{graphicx} 
  \DeclareGraphicsExtensions{.eps}
  
  \usepackage{tikz}
  \usepackage{pgfplots,pgfplotstable}
  \usetikzlibrary{pgfplots.groupplots}
  \pgfplotsset{compat=1.5}

\usepackage{multirow}

\usepackage{amsmath}

\usepackage{textcomp}
\usetikzlibrary{shapes,arrows}

\usepackage{empheq} 

\usepackage{stfloats} 

\usepackage[hyphens]{url}
\usepackage[hidelinks]{hyperref}
\hypersetup{colorlinks=true,linkcolor=black,citecolor=black,filecolor=black,urlcolor=black,breaklinks=true}
\usepackage{cite} 

\usepackage{color}

\usepackage{cancel} 

\usepackage{array} 

\usepackage[export]{adjustbox} 

\usepackage[normalem]{ulem} 

\ifCLASSINFOpdf
\else
\fi
\makeatletter
\def\bstctlcite{\@ifnextchar[{\@bstctlcite}{\@bstctlcite[@auxout]}}
\def\@bstctlcite[#1]#2{\@bsphack
 \@for\@citeb:=#2\do{%
   \edef\@citeb{\expandafter\@firstofone\@citeb}%
   \if@filesw\immediate\write\csname #1\endcsname{\string\citation{\@citeb}}\fi}%
 \@esphack}
\makeatother

\begin{document}
%

\bstctlcite{IEEEexample:BSTcontrol}

\title{Optimal Portfolio of Distinct Frequency Response\\Services in Low-Inertia Systems}
%
%
%

\author{Luis~Badesa,~\IEEEmembership{Student Member,~IEEE,}
        Fei~Teng,~\IEEEmembership{Member,~IEEE,}
        and~Goran~Strbac,~\IEEEmembership{Member,~IEEE}
\thanks{A portion of this work has been supported by project IDLES under grant EP/R045518/1.}
\thanks{
The authors are with the Department of Electrical and Electronic Engineering, Imperial College London, SW7 2AZ London, U.K. (email: luis.badesa@imperial.ac.uk, f.teng@imperial.ac.uk, g.strbac@imperial.ac.uk).}
%
}

%
%

\markboth{IEEE Transactions on Power Systems, May~2020}%
{Shell \MakeLowercase{\textit{et al.}}: Bare Demo of IEEEtran.cls for IEEE Journals}
%



\maketitle

\begin{abstract}
A reduced level of system inertia due to renewable integration increases the need for cost-effective provision of ancillary services, such as Frequency Response (FR). In this paper a closed-form solution to the differential equation describing frequency dynamics is proposed, which allows to obtain frequency-security algebraic constraints to be implemented in optimization routines. This is done while considering any finite number of FR services with distinguished characteristics, such as different delivery times and activation delays. The problem defined by these frequency-security constraints can be formulated as a Mixed-Integer Second-Order Cone Program (MISOCP), which can be efficiently handled by off-the-shelf conic optimization solvers. This paper also takes into account the uncertainty in inertia contribution from the demand side by formulating the frequency-security conditions as chance constraints, for which an exact convex reformulation is provided. Finally, case studies highlighting the effectiveness of this frequency-secured formulation are presented.
\end{abstract}


%
\IEEEpeerreviewmaketitle

\section*{Nomenclature}
\addcontentsline{toc}{section}{Nomenclature}

\begin{IEEEdescription}[\IEEEusemathlabelsep\IEEEsetlabelwidth{$\textrm{RoCoF}\textrm{ma}$}]
\setlength\itemsep{0.5em}
\item[$\alpha$] Probability of meeting the nadir constraint.
\item[$\Delta f(t)$] Time-evolution of post-fault frequency deviation from nominal state.
\item[$\Delta f_{\textrm{max}}$] Maximum admissible frequency deviation at the nadir (Hz).
\item[$\eta$] Probability of meeting the RoCoF constraint.
\item[$\Phi(\cdot)$] Standard normal cumulative distribution function.
\item[$\sigma$] Standard deviation of inertia from the demand side (MW$\cdot$s).
\item[$\tau_s$] Time-constant of FR provider $s$ (s).
\item[BESS] Battery Energy Storage Systems.
\item[CCGT] Combined Cycle Gas Turbine.
\item[$\textrm{D}$] Load damping (MW/Hz).
\item[EFR] Enhanced Frequency Response.
\item[$f_0$] Nominal frequency of the power grid (Hz).
\item[FR] Frequency Response.
\item[$\textrm{FR}(t)$] Time-evolution of aggregated system FR (MW).
\item[$g(H_\textrm{D})$] Auxiliary function for random variable $H_\textrm{D}$.
\item[GB] Great Britain.
\item[$H$] Decision variable, system inertia from thermal generators (MW$\cdot \textrm{s}$).
\item[${H_\textrm{D}}$] Random variable, system inertia from the demand side (MW$\cdot \textrm{s}$).
\item[$\textrm{H}_\mu$] Forecast for inertia from demand (MW$\cdot$s).
\item[$i,\,\, j,\,\, n$] All-purpose indices.
\item[$k,\,\, \mathcal{K}$] Index, Set of FR services fully delivered by the frequency nadir.
\item[$\textrm{K}_s$] Droop gain for FR provider $s$ (Hz/MW).
\item[$l,\,\, \mathcal{L}$] Index, Set of FR services ramping up by the frequency nadir.
\item[MISOCP] Mixed-Integer Second-Order Cone Program.
\item[$\mathcal{N}(\cdot,\cdot)$] Normal distribution.
\item[OCGT] Open Cycle Gas Turbine.
\item[$\textrm{\textbf{P}}(\cdot)$] Probability operator.
\item[$P_{\textrm{L}}$] Decision variable, largest power infeed (MW).
\item[$\textrm{P}_{\textrm{L}}^{\textrm{max}}$] Upper bound for $P_{\textrm{L}}$ (MW).
\item[PFR] Primary Frequency Response.
\item[$R_s$] Decision variable, maximum FR provision from service $s$ (MW).
\item[RoCoF] Rate-of-Change-of-Frequency.
\item[$\textrm{RoCoF}_{\textrm{max}}$] Maximum admissible RoCoF (Hz/s).
\item[$s,\,\, \mathcal{S}$] Index, Set of all FR services.
\item[$\left|\mathcal{S}\right|$] Cardinality of set $\mathcal{S}$.
\item[SOC] Second-Order Cone.
\item[SUC] Stochastic Unit Commitment.
\item[$\textrm{T}_s$] Delivery time of FR service $s$ (s).
\item[$\textrm{T}_{\textrm{del},s}$] Delay in provision of FR service $s$ (s).
\item[UC] Unit Commitment.
\end{IEEEdescription}

\newpage
\section{Introduction}
%
%
%
%
\IEEEPARstart{I}{ncreasing} penetration of renewable energy in power grids introduces many challenges, such as uncertainty and variability in generation. Furthermore, most renewable sources, including wind and photovoltaic, do not contribute to system inertia due to being decoupled from the grid by power electronic converters. While system inertia and Frequency Response (FR) were services widely available in grids dominated by thermal generation as by-products of energy production, the increasing scarcity of these services in low-carbon system increments the costs associated to their provision \cite{FeiISGT2017}. These frequency services are necessary to contain the frequency drop after a generation outage, in order to avoid the tripping of Rate-of-Change-of-Frequency (RoCoF) relays and/or the activation of Under-Frequency Load Shedding.  

In this context, several works have studied optimal strategies to provide inertia and FR \cite{ElaI,ERCOT_EFR,EDF_cuts,OPFChavez,IowaThesis,FeiStochastic,LinearizedUC,PricingElaZhang,VincenzoEFR,LuisEFR,UCFaroe}. Authors in \cite{ElaI} studied the design of an ancillary service market for FR, based on solving a constrained optimization problem that guarantees frequency security. The frequency-security constraints were obtained heuristically from dynamic simulations of the system. A similar approach was used in \cite{ERCOT_EFR,EDF_cuts}, while the other works focus on deducing the frequency-security region by analytically solving the differential equation describing post-fault frequency dynamics. 

To tackle the declining system inertia, FR services with faster delivery have been introduced by system operators \cite{NationalGridPLossInertia}. However, a fundamental question yet to be answered is how to optimize the portfolio of FR services from providers with diverse characteristics under different system conditions. References \cite{OPFChavez,IowaThesis,FeiStochastic,LinearizedUC,PricingElaZhang} aggregate the response from all FR providers uniformly, therefore only allowing to consider a single FR service. References \cite{LuisEFR,VincenzoEFR} co-optimize the two FR services defined in Great Britain (GB) up to date, namely Enhanced Frequency Response (EFR) delivered by one second after the outage, and Primary Frequency Response (PFR) delivered ten seconds after the outage. The only work that integrates different dynamics from generic FR providers is \cite{UCFaroe}, by considering an affine conservative approximation of the frequency dynamics until reaching the nadir. Furthermore, the impact of the activation delays in FR
has been shown to affect the system frequency stability \cite{QitengDelays}, which has not been considered in any of the above works. 

In the available literature, system inertia has been assumed either fixed in the optimization of FR \cite{ERCOT_EFR,OPFChavez} or fully controllable by the Unit Commitment (UC) \cite{FeiStochastic,LuisEFR}. In fact, a considerable amount of inertia contribution is available from demand, which is not controllable and can only be forecasted: assuming a central authority clearing a pool market using a frequency-secured UC, the inertia from demand can only be forecasted with certain accuracy. Given the risk aversion of system operators, it is necessary to explicitly model such uncertainty, in addition to the uncertainty associated with renewable generation.
Furthermore, optimally scheduling the largest online unit has been demonstrated to provide both economic and emission savings \cite{OMalleyDeload,LuisPESGM2018}, which also needs to be co-optimized along with other frequency services.

Given this background, this paper develops an efficient frequency-constrained optimization framework that recognizes and appropriately values the different dynamics of FR services. The contributions are three-fold: 
\begin{enumerate}
	\item  This paper proposes closed-form conditions to optimize, for the first time, any finite number of FR services with diverse dynamics while considering any combination of activation delays. The proposed frequency constraints are formulated as a Mixed-Integer Second-Order Cone Program (MISOCP), which can be efficiently solved by taking advantage of the recent development of conic optimization software.
    \item The uncertainty associated with the inertia contribution from the demand side is explicitly modelled in the optimization problem through chance constraints. A convex-reformulation of the chance constraints allows to maintain the problem as an MISOCP. 
    \item The frequency-constrained Stochastic Unit Commitment (SUC) model is applied to several case studies, which highlight the benefits of a frequency-security framework allowing to co-optimize a diverse portfolio of services. 
\end{enumerate}

The rest of this paper is organized as follows: Section \ref{SectionFrequency} provides the deduction of frequency-security constraints. The proposed analytical model for guaranteeing frequency security is validated through dynamic simulations in Section \ref{SectionValidation}. Section \ref{SectionCaseStudies} includes the results from several relevant case studies, while Section \ref{SectionConclusion} gives the conclusion and proposes future lines of work.

\section{Closed-Form Conditions for Secure Post-Fault Frequency Dynamics} \label{SectionFrequency}

In this section the closed-form conditions for frequency security are deduced, that allow to map the sub-second dynamics of transient frequency to any desired timescale, such as the typical minutes to hour resolution of a UC. Frequency security is respected if sufficient inertia and FR services are available at the time of the generation outage to contain and recover the system frequency. 

The conditions for a secure post-fault frequency evolution can be deduced from solving the swing equation \cite{KundurBook}:
\begin{equation} \label{SwingEqMultiFR}
\frac{2(H+H_\textrm{D})}{f_0}\frac{\textrm{d} \Delta f(t)}{\textrm{d} t} = \textrm{FR}(t) - P_{\textrm{L}}
\end{equation}
Eq. (\ref{SwingEqMultiFR}) assumes the loss of the largest power infeed, therefore representing the common $N-1$ reliability requirement in power systems. Load damping has been neglected, as the damping level will be significantly reduced in future power systems that are increasingly dominated by power electronics \cite{OPFChavez}. Note that eq. (\ref{SwingEqMultiFR}) considers the uniform-frequency model, which assumes frequency to be approximately equal in all buses of the grid; electric frequency is therefore considered as a system-wide magnitude, while the electric network could be simultaneously considered for any other magnitude such as nodal voltages.

The system inertia is aggregated from the inertia provided by all devices, including thermal generators and certain loads, and it includes two components: $H$ is the controllable term in the UC optimization, as it is a decision variable which depends on the generators that are scheduled to be online; on the other hand, $H_\textrm{D}$ is the inertia contribution from demand, which is assumed here to be non-controllable but can be forecasted for time-periods in the near future. The largest possible generation outage $P_{\textrm{L}}$ can be considered as a decision variable, which would involve part-loading a large generating unit or interconnector \cite{LuisPESGM2018,LuisEFR}. 

Function $\textrm{FR}(t)$ in (\ref{SwingEqMultiFR}) represents the frequency control for Frequency Response, a power injection from several system devices following the outage. This function is modelled in (\ref{definitionMultiFR}) to consider $|\mathcal{S}|$ different FR services, in which each FR service is associated with a characteristic delivery time, therefore allowing to recognize diverse dynamics in FR delivery from different providers.

\begin{subequations} \label{definitionMultiFR}
\begin{empheq}[left ={\hspace{-5pt}\textrm{FR}(t)=\empheqlbrace}]{alignat=3}
 & \sum_{s\in\mathcal{S}}\frac{R_s}{\textrm{T}_{s}}\cdot t && \mbox{if $t\leq\textrm{T}_1$} \tag{\ref{definitionMultiFR}.1}\\[-2pt]
 & R_1 + \sum_{s=2}^{|\mathcal{S}|}\frac{R_s}{\textrm{T}_{s}}\cdot t && \mbox{if $\textrm{T}_1<t\leq \textrm{T}_2$} \tag{\ref{definitionMultiFR}.2}\\[-3pt]
 & & \mathrel{\cdots} & & \nonumber \\[-3pt]
  & \sum_{s=1}^{|\mathcal{S}|-1}R_s + \frac{R_{|\mathcal{S}|}}{\textrm{T}_{|\mathcal{S}|}} t && \mbox{if $\textrm{T}_{|\mathcal{S}|-1}<t\leq \textrm{T}_{|\mathcal{S}|}$} {\hspace{-4mm}\tag{\ref{definitionMultiFR}.$|\mathcal{S}|$}} \\[1pt]
 & \sum_{s\in\mathcal{S}}R_s && \mbox{if $t> \textrm{T}_{|\mathcal{S}|}$} \tag{\ref{definitionMultiFR}.($|\mathcal{S}|$+1)}
\end{empheq}
\end{subequations}

The delivery time $\textrm{T}_s$ of a service $s$ is the time by which full FR capacity for the service is delivered. Piecewise function (\ref{definitionMultiFR}) models the delivery of FR from each provider $s$ as ramping up during the interval $t \in (0,\textrm{T}_{s}]$, and constant for $t > \textrm{T}_{s}$. This approach, proposed in \cite{OPFChavez} for a single FR service, is guaranteed to conservatively approximate any controller proportional to frequency deviation, as is further demonstrated in Section \ref{SectionValidation}. Note that FR services are ordered from service $1$ up to service $|\mathcal{S}|$ in increasing delivery time, i.e. service $\textrm{FR}_1$ is the fastest and service $\textrm{FR}_{|\mathcal{S}|}$ is the slowest.

\subsection{Deducing Frequency-Security Constraints} \label{SectionSolveSwing}

By solving (\ref{SwingEqMultiFR}), the three constraints that guarantee a secure post-fault frequency evolution can be obtained. The RoCoF constraint is deduced following the standard from National Grid \cite{NGrocofDefinition}, which established that the level of system inertia must be sufficient to limit the highest instantaneous RoCoF at $t=0$:
\begin{equation} \label{RocofConstraint}
\left|\mbox{RoCoF}(t=0)\right| = \frac{P_{\textrm{L}}\cdot f_0}{2 (H+H_\textrm{D})} \leq \mbox{RoCoF}_{\textrm{max}}
\end{equation}

For frequency to stabilize eventually after the fault, the amount of FR available must be at least equal to the generation outage. In other words, the steady-state constraint is obtained by setting RoCoF to zero in (\ref{SwingEqMultiFR}) and considering that every FR service has been fully delivered: 
\begin{equation} \label{qssMultiFR}
\sum_{s\in\mathcal{S}}R_s \geq P_{\textrm{L}}
\end{equation}

Since $\textrm{FR}(t)$ defined in (\ref{definitionMultiFR}) is a piecewise function, the nadir constraint depends on the time-interval when the nadir occurs. For the nadir to take place in a given time-interval $t \in \left[\textrm{T}_{n-1},\textrm{T}_n\right)$, the two following conditions must be met:
\begin{multline} \label{ConditionNadirNoD}
\quad\Biggl(\sum_{i=1}^{n-1}R_i + \sum_{j=n}^{|\mathcal{S}|} R_j \frac{\textrm{T}_{n-1}}{\textrm{T}_{j}} \leq P_{\textrm{L}}\Biggr) \quad \textrm{and} \\[-3pt]
\Biggl(\sum_{i=1}^n R_i + \sum_{j=n+1}^{|\mathcal{S}|} R_j \frac{\textrm{T}_{n}}{\textrm{T}_{j}} > P_{\textrm{L}}\Biggr)
\end{multline}
Condition (\ref{ConditionNadirNoD}) states that the power injected from FR becomes greater than the power loss $P_\textrm{L}$ not before $\textrm{T}_{n-1}$ and no later than $\textrm{T}_n$. Before $\textrm{T}_{n-1}$, only the fastest $n-1$ FR services have been fully delivered, while the rest are still ramping up; after $\textrm{T}_n$, the $n^\textrm{th}$ service has been fully delivered as well.

The solution of (\ref{SwingEqMultiFR}) for $t \in \left[\textrm{T}_{n-1},\textrm{T}_n\right)$ is:
\begin{equation} \label{Sol_Delta_f}
\Delta \hspace{-0.5mm}f(t) \hspace{-0.5mm}=\hspace{-0.5mm} \frac{f_0}{2 (H\hspace{-1mm}+\hspace{-1mm}H_\textrm{D})}\hspace{-1mm}\left[\sum_{j=n}^{|\mathcal{S}|}\frac{R_j}{2\textrm{T}_j} t^2 \hspace{-1mm}+\hspace{-1mm} \sum_{i=1}^{n-1} \hspace{-0.5mm}R_i \biggl(\hspace{-0.6mm}t\hspace{-0.5mm}-\hspace{-0.5mm}\frac{\textrm{T}_i}{2}\hspace{-0.3mm}\biggr) \hspace{-0.7mm}-\hspace{-0.5mm} P_{\textrm{L}} \hspace{-1mm}\cdot\hspace{-0.5mm} t\right]
\end{equation}
The time within $t \in \left[\textrm{T}_{n-1},\textrm{T}_n\right)$ at which nadir is exactly reached is given by setting RoCoF to zero in (\ref{SwingEqMultiFR}) for that given time interval:
\begin{equation} \label{t_nadir}
t_\textrm{nadir} = \frac{P_{\textrm{L}} - \sum_{i=1}^{n-1}R_i}{\sum_{j=n}^{|\mathcal{S}|}R_j/\textrm{T}_{j}} 
\end{equation}
By substituting (\ref{t_nadir}) into (\ref{Sol_Delta_f}), the condition for respecting the nadir requirement can be deduced:
\begin{equation} \label{nadir_req}
\left| \Delta f_{\textrm{nadir}} \right| = \left| \Delta f(t=t_\textrm{nadir})  \right| \leq \Delta f_{\textrm{max}}
\end{equation}
Finally, expanding the expression in (\ref{nadir_req}) and enforcing the conditions for $t_\textrm{nadir}$ to occur during time-interval $t \in \left[\textrm{T}_{n-1},\textrm{T}_n\right)$, the nadir constraint is obtained as: 
\vspace{4mm}

if $\; \left(\sum_{i=1}^{n-1}R_i + \sum_{j=n}^{|\mathcal{S}|} R_j \frac{\textrm{T}_{n-1}}{\textrm{T}_{j}} \leq P_{\textrm{L}}\right)$ and \\
\hspace*{10mm}$\left(\sum_{i=1}^n R_i + \sum_{j=n+1}^{|\mathcal{S}|} R_j \frac{\textrm{T}_{n}}{\textrm{T}_{j}} > P_{\textrm{L}}\right)$ then enforce:
\begin{equation} \label{ConditionalSOCconstraints}
\Biggl(\underbrace{\vphantom{\sum_{j=n}^S} \frac{H\hspace{-1mm}+\hspace{-1mm}H_\textrm{D}}{f_0} - \sum_{i=1}^{n-1} \frac{R_i \textrm{T}_{i}}{4 \Delta f_{\textrm{max}}}}_\text{$= x_1$}\Biggr) \underbrace{\sum_{j=n}^{|\mathcal{S}|} \frac{R_j}{\textrm{T}_{j}}}_\text{$= x_2$}
\geq 
\frac{(P_{\textrm{L}}\hspace{-1mm}-\sum_{i=1}^{n-1} R_i)^2}{4 \Delta f_{\textrm{max}}}
\end{equation}

As one must consider the possibility of nadir occurring at any time $t \in \left[0,\textrm{T}_{|\mathcal{S}|}\right)$ (note that the nadir must occur before $\textrm{T}_{|\mathcal{S}|}$, as otherwise the steady-state constraint (\ref{qssMultiFR}) would not hold), $|\mathcal{S}|$ different nadir constraints must be defined, corresponding to each time-interval $\left[\textrm{T}_\textrm{s-1},\textrm{T}_\textrm{s}\right) \; \forall s \in \mathcal{S}$. Only one constraint will be enforced, which is the constraint for which the if-statement in (\ref{ConditionalSOCconstraints}) is met. Note that conditional statements in optimization can be implemented using a big-M formulation with auxiliary binary decision variables \cite{ConejoOptimizationBook}. 

Constraints (\ref{ConditionalSOCconstraints}) are non-linear but are in fact rotated Second-Order Cones (SOCs), therefore convex constraints as  $x_1$ and $x_2$ in (\ref{ConditionalSOCconstraints}) are non-negative. SOC Programming generalizes Linear Programming, and recently developed interior-point methods allow to efficiently solve these types of conic optimization problems to global optimality \cite{BoydConvex}. Furthermore, SOC Programs are the highest class of conic problems whose mixed-integer counterpart can be solved to global optimality using commercial optimization packages. Since the nadir constraints in (\ref{ConditionalSOCconstraints}) introduce binary variables for implementing the conditional statements, the resulting optimization problem is an MISOCP.

In conclusion, constraints (\ref{RocofConstraint}), (\ref{qssMultiFR}) and (\ref{ConditionalSOCconstraints}) guarantee frequency security in a power system, while considering the dynamics of any finite number $|\mathcal{S}|$ of different FR providers.

\subsection{Considering Activation Delays in Certain FR Services} \label{SectionDelays}

\begin{figure}
\raggedright
\hspace*{-0.45cm}
    \input{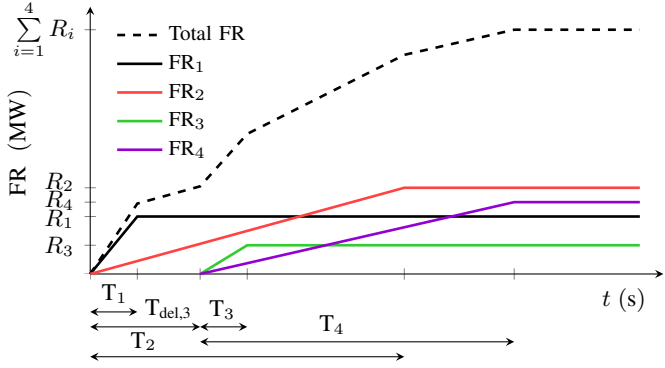}
    \caption{Time evolution of four distinct FR services: the first two start ramping up at the very moment of the generation outage and the other two have an activation delay. Note that $\textrm{T}_{\textrm{del},4}=\textrm{T}_{\textrm{del},3}$ in this case. The total system FR(t) as defined in eq. (\ref{SwingEqMultiFR}) is given by the dashed line.}
	\label{ExampleStaticDynamic}
\end{figure}

In Section \ref{SectionSolveSwing}, every FR service is considered in (\ref{definitionMultiFR}) to start ramping up at the very moment of the generation outage. Those FR services would therefore react to any deviation from nominal frequency in the grid, not necessarily caused by the loss of a large power infeed. Here the model is generalized to account for some FR services which start providing FR some time after the outage:
\begin{subequations} \label{definitionStaticFR}
\begin{empheq}[left ={\hspace{-5pt}\textrm{FR}_{i}(t)=\empheqlbrace}]{alignat=3}
 & 0 &\, &\quad && \mbox{if $t\leq\textrm{T}_{\textrm{del},i}$} \tag{\ref{definitionStaticFR}.1}\\[-2pt]
 & \frac{R_i}{\textrm{T}_{i}} (t-\textrm{T}_{\textrm{del},i}) &\, &\quad && \mbox{if $\textrm{T}_{\textrm{del},i}<t\leq \textrm{T}_i$} \tag{\ref{definitionStaticFR}.2}\\[-2pt]
 & R_i &\, &\quad && \mbox{if $t > \textrm{T}_i$} \tag{\ref{definitionStaticFR}.3}
\end{empheq}
\end{subequations}
The FR service defined in (\ref{definitionStaticFR}) is activated $\textrm{T}_{\textrm{del},i}$ seconds after the fault, a delay which can be driven by either the frequency deadband of a droop control or the communication delay of an activation signal sent to the FR provider. For the following deductions in this section, $\textrm{FR}(t)$ as defined in (\ref{definitionMultiFR}) may now include FR services with an activation delay as the one defined in (\ref{definitionStaticFR}). An example considering four FR services is included in Fig. \ref{ExampleStaticDynamic}, where services $\textrm{FR}_1$ and $\textrm{FR}_2$ start ramping up at the very moment of the fault (the fault is assumed to happen at $t=0$) while services $\textrm{FR}_3$ and $\textrm{FR}_4$ have an activation delay.

The RoCoF and steady-state constraints, (\ref{RocofConstraint}) and (\ref{qssMultiFR}), remain unchanged while the nadir constraints must be updated if some FR services have an activation delay. Following the same procedure as in Section \ref{SectionSolveSwing}, the swing equation is solved for the different time-intervals, yielding the following nadir constraint:
\begin{multline} \label{SOCnadirStatic}
\hspace{-4mm}\Biggl( \frac{H\hspace{-1mm}+\hspace{-1mm}H_\textrm{D}}{f_0} \underbrace{ - \sum_{k \in \mathcal{K}} \frac{R_k (\textrm{T}_{k}\hspace{-1mm}+\hspace{-1mm}2\textrm{T}_{\textrm{del},k})}{4 \Delta f_{\textrm{max}}}  + \sum_{l \in \mathcal{L}} \frac{R_l \textrm{T}_{\textrm{del},l}^2/\textrm{T}_{l}}{4 \Delta f_{\textrm{max}}} }_\text{$= y_1$}\Biggr) \underbrace{\sum_{l \in \mathcal{L}} \frac{R_l}{\textrm{T}_{l}}}_\text{$= y_2$} \\
\geq \underbrace{\frac{(P_{\textrm{L}} - \sum_{k \in \mathcal{K}} R_k + \sum_{l \in \mathcal{L}} R_l \textrm{T}_{\textrm{del},l}/\textrm{T}_{l})^2}{4 \Delta f_{\textrm{max}}}}_\text{$= y_3^2$}
\end{multline}
Note that if every FR service starts ramping up exactly when the fault occurs, i.e. $\textrm{T}_{\textrm{del},k}=\textrm{T}_{\textrm{del},l}=0 \;\; \forall k\in\mathcal{K},\forall l\in\mathcal{L}$, (\ref{SOCnadirStatic}) reduces to (\ref{ConditionalSOCconstraints}). Therefore, (\ref{SOCnadirStatic}) generalizes (\ref{ConditionalSOCconstraints}) allowing to consider any combination of activation delays for FR services. The problem defined is still an MISOCP, since (\ref{SOCnadirStatic}) is a rotated SOC. In a similar fashion as in (\ref{ConditionalSOCconstraints}), as many nadir constraints as intervals defined by the piecewise $\textrm{FR}(t)$ must be included, along with the corresponding conditional statements for nadir to occur in that interval. For the example in Fig. \ref{ExampleStaticDynamic}, each time-interval is delimited by a tick in the x-axis.

\subsection{Uncertainty in Inertia Contribution from Demand} \label{SectionChance}

In this section the use of chance constraints is proposed, that allow to take into account the inertia contribution from demand subject to forecasting errors.
The inertia from demand, $H_\textrm{D}$, is considered as a random variable for which a forecast is available, along with a distribution on the forecasting error. In order to account for this uncertainty in the frequency-security conditions deduced in Sections \ref{SectionSolveSwing} and \ref{SectionDelays}, the RoCoF and nadir constraints are modified to become chance constraints, i.e. constraints that must be met above a pre-defined probability. Here an exact convex reformulation of the non-convex chance constraints is provided, to allow the system operator to limit the risk of violating each frequency constraint. The error in inertia forecasting is assumed to follow a Gaussian distribution (but any log-concave probability density function still makes the following deductions valid):
\begin{equation}
H_\textrm{D} \sim \mathcal{N} (\textrm{H}_\mu,\sigma^2)
\end{equation}

The chance constraint for meeting the RoCoF requirement is given by the following non-convex constraint, based on (\ref{RocofConstraint}):
\begin{equation} \label{RocofConstraintChance}
\textrm{\textbf{P}}\left( H_\textrm{D}  \geq \frac{P_{\textrm{L}}\cdot f_0}{2\cdot \mbox{RoCoF}_{\textrm{max}}} - H \right) \geq \eta
\end{equation}
Note that $H$ is a decision variable since it is the inertia contribution from the generators scheduled to be online in the UC, and therefore $H$ is not subject to uncertainty.

Since the constraint inside the probability operator in (\ref{RocofConstraintChance}) is linear, and making use of the log-concave property of the normal distribution for $H_\textrm{D}$ \cite{SaphiroChance}, the non-convex chance constraint (\ref{RocofConstraintChance}) is equivalent to:
\begin{equation} 
\Phi \left(\frac{-\frac{P_{\textrm{L}}\cdot f_0}{2\cdot \textrm{RoCoF}_{\textrm{max}}} + H + \textrm{H}_\mu}{\sigma} \right) \geq \eta
\end{equation}

Therefore, the exact linear reformulation of the RoCoF chance constraint (\ref{RocofConstraintChance}) is obtained making use of the inverse cumulative distribution function:
\begin{equation} \label{reformulationRocofChance1}
\frac{ -\frac{P_{\textrm{L}}\cdot f_0}{2\cdot \textrm{RoCoF}_{\textrm{max}}} + H + \textrm{H}_\mu}{\sigma} \geq \Phi^{-1}(\eta) 
\end{equation}
Rearranging (\ref{reformulationRocofChance1}) for clarity:
\begin{equation} \label{reformulationRocofChance}
H + \textrm{H}_\mu - \Phi^{-1}(\eta) \sigma \geq \frac{P_{\textrm{L}}\cdot f_0}{2\cdot \mbox{RoCoF}_{\textrm{max}}} 
\end{equation}

The chance constraint for the nadir requirement, using the notation in (\ref{SOCnadirStatic}), is:
\begin{equation} \label{NadirConstraintChance}
\textrm{\textbf{P}}\biggl[ \biggl(\frac{H+H_\textrm{D}}{f_0} + y_1\biggr) y_2 \geq y_3^2
\biggr] = \textrm{\textbf{P}}\left[ g(H_\textrm{D}) \leq 0 \right] \geq \alpha
\end{equation}
Where function $g(H_\textrm{D})$ is given by:
\begin{equation}
g(H_\textrm{D})=-\left(\frac{H+H_\textrm{D}}{f_0} + y_1 \right) y_2 + y_3^2
\end{equation}

Function $g(H_\textrm{D})$ is linear with respect to the random variable $H_\textrm{D}$, therefore it also follows a normal distribution \cite{RandomVariableBook}:
\begin{equation}
g(H_\textrm{D}) \sim \mathcal{N} \biggl(\underbrace{-\biggl[\frac{H+\textrm{H}_\mu}{f_0}+y_1 \biggr]y_2 + y_3^2}_\text{$=\mu_g$}, \,\, \underbrace{\sigma^2 \biggl[\frac{y_2}{f_0}\biggr]^2}_\text{$=\sigma_g^2$} \biggr)
\end{equation}

Again making use of the log-concave property of the normal distribution, (\ref{NadirConstraintChance}) becomes:
\begin{equation} \label{reformulationNadirChance1}
\Phi \left(\frac{0-\mu_g}{\sigma_g} \right) \geq \alpha 
\end{equation}
Expanding (\ref{reformulationNadirChance1}):
\begin{equation} \label{reformulationNadirChance2}
\biggl(\frac{H+\textrm{H}_\mu}{f_0}+y_1\biggr)y_2 - y_3^2 \geq \Phi^{-1}(\alpha)\sigma \frac{y_2}{f_0}
\end{equation}
Finally, rearranging (\ref{reformulationNadirChance2}):
\begin{equation} \label{reformulationNadirChance}
\biggl(\frac{H+\textrm{H}_\mu-\Phi^{-1}(\alpha)\sigma}{f_0}+y_1\biggr)y_2 \geq y_3^2 
\end{equation}

Constraint (\ref{reformulationNadirChance}) is a rotated SOC, therefore it provides a convex reformulation of the chance constraint (\ref{NadirConstraintChance}), using the notation for the linear expressions $y_1$, $y_2$ and $y_3$ from (\ref{SOCnadirStatic}). As in (\ref{SOCnadirStatic}), any combination of distinct FR services with or without activation delays can be considered.

\section{Validation and Applicability of the Frequency-Security Constraints} \label{SectionValidation}

The frequency-security constraints obtained in Section \ref{SectionFrequency} are purely mathematical deductions from the swing equation (\ref{SwingEqMultiFR}), and hence guaranteed to provide the security region entailing no approximation. However, the assumptions for function $\textrm{FR}(t)$ in (\ref{definitionMultiFR}) are conservative, as considering detailed frequency controls in the swing equation would impede to solve it algebraically, and therefore no closed-form frequency-security conditions could be obtained. In this section it is demonstrated that the assumptions for $\textrm{FR}(t)$ do indeed underestimate an actual frequency droop control, a demonstration based on comparing $\textrm{FR}(t)$ with an actual dynamic simulation to which the solution of a frequency-constrained optimization is fed. 

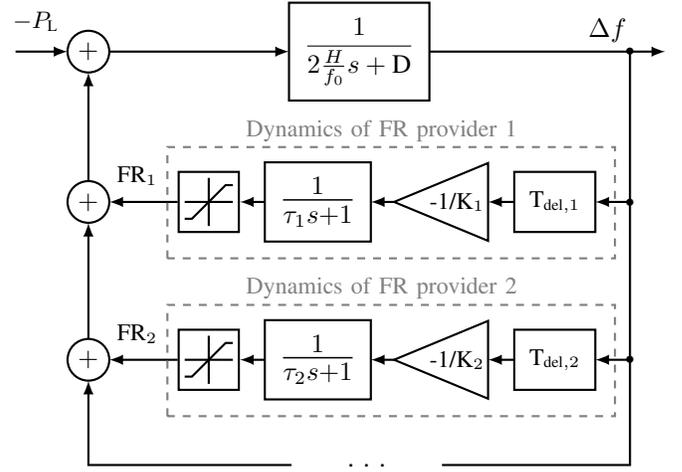
\begin{figure} 
\raggedright

\tikzset{%
  block/.style    = {draw, thick, rectangle, minimum height = 3em,
    minimum width = 3em, inner sep=2mm, outer sep=0mm,},
  gain/.style     = {draw, thick, isosceles triangle, minimum height = 0.5em, inner sep=0.5mm, outer sep=0mm,
     isosceles triangle apex angle=45, shape border rotate=-180},
  sum/.style      = {draw, circle, inner sep=2pt, outer sep=0pt, node distance = 1cm}, 
  input/.style    = {coordinate}, 
  output/.style   = {coordinate}, 
  saturation/.style={%
    draw, 
    path picture={
      \pgfpointdiff{\pgfpointanchor{path picture bounding box}{north east}}%
        {\pgfpointanchor{path picture bounding box}{south west}}
      \pgfgetlastxy\x\y
      \tikzset{x=\x*.4, y=\y*.4}
      %
      \draw (-1,0) -- (1,0) (0,-1) -- (0,1); 
      \draw (-1,-.7) -- (-.7,-.7) -- (.7,.7) -- (1,.7);
    }
  }
}
\newcommand{\suma}{$+$}

\begin{tikzpicture}[auto, thick, node distance=1cm, >=latex]

\draw
    node [] (input1) {}
    node at (1.1,0)[sum] (sum1) {\suma}
    node at (4.7,0)[block] (InertiaBlock) {$\dfrac{1}{2\frac{H}{f_0} s + \textrm{D} }$}
    node at (8.9,0)[] (output1) {}
    node at (8.3,0.1)[](derivation1) {}
    
    node at (7.3,-2)[block, minimum height=0.8cm] (delay1) {\small $\textrm{T}_{\textrm{del},1}$}
    node at (6,-2)[gain] (droop1) {\small -1/$\textrm{K}_1$}
    node at (4.15,-2)[block] (governor1) {\Large $\frac{1}{\tau_1 s + 1}$}
    node at (2.7,-2)[saturation, minimum width=0.8cm, minimum height=0.8cm] (sat1) {}
    node at (8.3,-1.9)[](derivation2) {}
    
    node at (7.3,-4.1)[block, minimum height=0.8cm] (delay2) {\small $\textrm{T}_{\textrm{del},2}$}
    node at (6,-4.1)[gain] (droop2) {\small -1/$\textrm{K}_2$}
    node at (4.15,-4.1)[block] (governor2) {\Large $\frac{1}{\tau_2 s + 1}$}
    node at (2.7,-4.1)[saturation, minimum width=0.8cm, minimum height=0.8cm] (sat2) {}

    node at (1.1,-2)[sum](sum2) {\suma}
    node at (8.3,-4)[](derivation3) {}
    node at (4.8,-5.5)[block, draw=none, minimum width=2cm] (etc) {. . .}
    node at (1.1,-4.1)[sum](sum3) {\suma}
    ;
    
    \draw[->](input1) -- node {} (sum1);
	\draw[->](sum1) -- node {} (InertiaBlock);
    \draw[->](InertiaBlock) -- node {\hspace{16mm}$\Delta f$} (output1);
    \draw[->](derivation1) |- node[near end]{} (delay1);
    \draw[->](delay1) -- node[near end]{} (droop1);
    \draw[->](droop1) -- node[near end]{} (governor1);
	\draw[->](governor1) -- node[near end]{} (sat1);
    \draw[->](sat1) -- node[near end]{} (sum2);
    \draw[->](derivation2) |- node[near end]{} (delay2);
    \draw[->](delay2) -- node[near end]{} (droop2);
    \draw[->](droop2) -- node[near end]{} (governor2);
	\draw[->](governor2) -- node[near end]{} (sat2);
    \draw[->](sat2) -- node[near end]{} (sum3);
    \draw[->](sum2) -- node[near end]{} (sum1);
    \draw[->](sum3) -- node[near end]{} (sum2);
    \draw[-](derivation3) |- node[near end]{} (etc);
    \draw[->](etc) -| node[near end]{} (sum3);

\draw
	node at (8.3,-0.035){\textbullet}
    node at (8.3,-2.035){\textbullet}
    node at (8.3,-4.135){\textbullet}
    ;
    
    \draw node at (0.4,0.4) {\small $-P_\textrm{L}$};
    
    \draw node at (1.75,-1.65) {\small $\textrm{FR}_1$};
    \draw node at (1.75,-3.75) {\small $\textrm{FR}_2$};
    
	\draw [color=gray, dashed,thick](2.15,-2.75) rectangle (8.1,-1.27);
	\node [color=gray] at (5,-1.05) {\small Dynamics of FR provider 1};
    \draw [color=gray, dashed,thick](2.15,-4.85) rectangle (8.1,-3.37);
	\node [color=gray] at (5,-3.15) {\small Dynamics of FR provider 2};
	
\end{tikzpicture}
\vspace{-20pt}
\caption{Block diagram for the simulation of the system frequency dynamics.}
\label{FigBlockDiagram}
\end{figure}

For the validation of the dynamic model for post-fault frequency, a generic case including four FR services has been considered: two FR services with no activation delay, $\textrm{FR}_1$ and $\textrm{FR}_2$, with delivery times $\textrm{T}_1=3\textrm{s}$ and $\textrm{T}_2=10\textrm{s}$; an FR service, $\textrm{FR}_3$, with $\textrm{T}_{\textrm{del},3}=0.5\textrm{s}$ and $\textrm{T}_3=5\textrm{s}$; and another FR service, $\textrm{FR}_4$, with $\textrm{T}_{\textrm{del},4}=1\textrm{s}$ and $\textrm{T}_4=8\textrm{s}$. An operating point exactly meeting the nadir constraint (\ref{SOCnadirStatic}) was fed into a dynamic simulation in MATLAB/Simulink, for which the dynamics of FR providers are modelled as in Fig. \ref{FigBlockDiagram}. The operating condition exactly meeting the nadir is $P_\textrm{L}=1.8\textrm{GW}$, $H=180\textrm{GWs}$, $R_1=0.2\textrm{GW}$, $R_2=0.98\textrm{GW}$, $R_3=0.5\textrm{GW}$, $R_4=0.6\textrm{GW}$. A load damping term of $\textrm{D}=0.15\textrm{GW/Hz}$ was added to the dynamic simulation, in order to analyze the impact of neglecting such support in the frequency constraints. 

\begin{figure} [!t]
\hspace*{0.02cm}
%
%
\definecolor{mycolor1}{rgb}{1.00000,0.26275,0.26275}%
\begin{tikzpicture}

\begin{axis}[%
axis lines = left,
width=2.65in,
height=1in,
at={(0in,0in)},
scale only axis,
y label style={at={(axis description cs:-0.2,0.25)},anchor=west},
xmin=0,
xmax=17,
xlabel style={font=\footnotesize},
xlabel={time (s)},
xtick={0,5,10,15},
ymin=-0.9,
ymax=0.05,
ylabel style={font=\footnotesize},
ylabel style={align=center},
ylabel={$\Delta f$ (Hz)},
axis background/.style={fill=white},
legend style={at={(0.03,0.97)}, anchor=north west, legend cell align=left, align=left, draw=white!15!black},legend style={font=\footnotesize},
legend style={fill opacity=0,text opacity=1,draw=none}, 
xticklabel style={font=\footnotesize},
yticklabel style={font=\footnotesize},
y axis line style = {-} 
]






\addplot [color=black, line width=1.0pt]
  table[row sep=crcr]{%
0	0\\
0.00142347427466447	-0.000355857447725201\\
0.00427042282399341	-0.00106754899569523\\
0.00711737137332235	-0.00177919169820546\\
0.0106891885417321	-0.00267195178713017\\
0.0254894827723628	-0.00637013738149229\\
0.0402897770029935	-0.0100663116317494\\
0.0581093459416577	-0.0145133186955386\\
0.0759289148803219	-0.0189564131829573\\
0.0937484838189861	-0.0233952419838344\\
0.115963627645627	-0.0289226256703495\\
0.138178771472267	-0.0344426000864307\\
0.167452065963889	-0.0417045704309435\\
0.196725360455511	-0.0489527066591083\\
0.237127374227608	-0.0589328921904214\\
0.277529387999705	-0.0688851921159406\\
0.339343188670331	-0.0840562124799395\\
0.417212072481069	-0.103068672351577\\
0.441265410238510	-0.108918730275395\\
0.465318747995951	-0.114757755702884\\
0.479125833089102	-0.118104441861739\\
0.492932918182253	-0.121447428419126\\
0.503415518963603	-0.123982998611409\\
0.513898119744953	-0.126516387081299\\
0.522161320623422	-0.128511830322499\\
0.530424521501891	-0.130505855107517\\
0.553458079619275	-0.136056394207900\\
0.576491637736659	-0.141595070585321\\
0.599525195854043	-0.147121290573137\\
0.622558753971427	-0.152634484800717\\
0.645592312088810	-0.158134107133313\\
0.674423962119507	-0.164998161350012\\
0.703255612150204	-0.171839154289927\\
0.740476520869819	-0.180635029982241\\
0.777697429589433	-0.189389045363703\\
0.824801031718298	-0.200404530580647\\
0.871904633847163	-0.211346859342797\\
0.930810629128707	-0.224923198041038\\
0.989716624410251	-0.238374981588111\\
1.00345987305324	-0.241495019159575\\
1.01720312169624	-0.244607929939352\\
1.02705321482014	-0.246834592167482\\
1.03690330794405	-0.249057513755280\\
1.06936791292857	-0.256356977421496\\
1.10183251791310	-0.263614261714173\\
1.13429712289762	-0.270828229527256\\
1.17517284206518	-0.279847888035933\\
1.21604856123274	-0.288795087337562\\
1.27280479451472	-0.301094334497608\\
1.32956102779670	-0.313245208520017\\
1.40565085278448	-0.329294097365907\\
1.48174067777226	-0.345058108416036\\
1.58070600350878	-0.365119699624341\\
1.67967132924531	-0.384666390411574\\
1.80403804942313	-0.408475746273475\\
1.92840476960094	-0.431422142981875\\
2.08078230151951	-0.458329052466964\\
2.23315983343808	-0.483879680482817\\
2.41663699103008	-0.512815600999922\\
2.60011414862209	-0.539733616267303\\
2.78359130621409	-0.564626463502907\\
3.01645310722537	-0.593315636031090\\
3.24931490823665	-0.618792977289408\\
3.53054517042089	-0.645376431376506\\
3.81177543260512	-0.667525153589593\\
4.09300569478935	-0.685433933109898\\
4.45172515261449	-0.702501168531355\\
4.81044461043962	-0.713564439149692\\
5.26920908357899	-0.719907879689913\\
5.72797355671835	-0.718599890740325\\
6.18673802985771	-0.710903649464679\\
6.64550250299707	-0.698060716652414\\
7.10426697613644	-0.681257519564554\\
7.56303144927580	-0.661598799329601\\
8.02179592241516	-0.640088814274155\\
8.48056039555453	-0.617618389212301\\
8.93932486869389	-0.594957812300528\\
9.39808934183325	-0.572754590422928\\
9.85685381497261	-0.551535297995011\\
10.3156182881120	-0.531710725397582\\
10.7743827612513	-0.513583583547109\\
11.2331472343907	-0.497358080816739\\
11.6919117075301	-0.483150754815681\\
12.1506761806694	-0.471002015359922\\
12.6094406538088	-0.460887930547617\\
13.0682051269482	-0.452731864148374\\
13.5269696000875	-0.446415646298775\\
13.9857340732269	-0.441790029675050\\
14.4444985463662	-0.438684248125592\\
14.9032630195056	-0.436914553508972\\
15.3620274926450	-0.436291658411036\\
15.8207919657843	-0.436627057310166\\
16.2795564389237	-0.437738236498721\\
};

\addplot [color=black, dashed]
table[row sep=crcr]{%
0 -0.8\\
23.9582891328490 -0.8\\
};

\end{axis}
\end{tikzpicture}%
    \caption{Post-fault frequency deviation from the dynamic simulation.}
	\label{FigFrequencySimulation}
\end{figure}
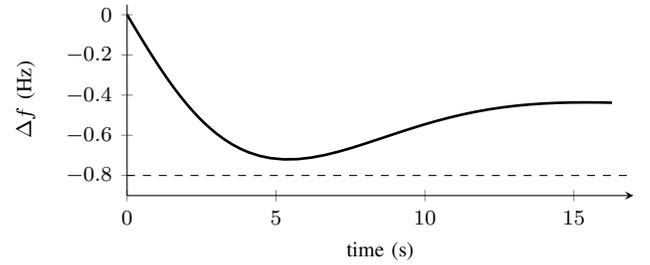

\begin{figure} [!t]
\hspace*{0.02cm}
    \input{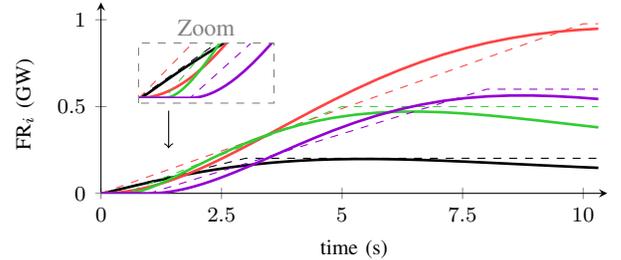}
    \caption{Time-evolution of FR obtained from the dynamic simulation considering the four different providers: $\textrm{FR}_1$ and $\textrm{FR}_2$ do not have an activation delay (black and red lines), while $\textrm{FR}_3$ and $\textrm{FR}_4$ have a delay (green and purple lines). The dashed lines represent the FR profile for each provider assumed in (\ref{definitionMultiFR}).}
	\label{FigFR_Simulation}
\end{figure}

The results of the dynamic simulation are shown in Fig. \ref{FigFrequencySimulation} and Fig. \ref{FigFR_Simulation}. The nadir is of $0.72\textrm{Hz}$, indeed above the $\Delta f_\textrm{max}=0.8\textrm{Hz}$ requirement. Although the nadir constraint was binding for this system condition, the $0.08\textrm{Hz}$ conservativeness in the simulation is due to neglecting the damping support and the linear-ramp assumption for FR delivery in (\ref{definitionMultiFR}). 

The validation of the frequency constraints has been performed for a particular system condition simply to illustrate the appropriateness of the proposed model. Even for contingencies of different sizes, this model guarantees post-fault frequency security as long as the ramps for FR are appropriately chosen to represent the frequency controls: the linear-ramp assumption for FR delivery in (\ref{definitionMultiFR}) can conservatively approximate a generic FR dynamics, such as the droop control in Fig. \ref{FigBlockDiagram}, as was demonstrated by \cite{OPFChavez}.

\subsection{Applicability to System Scheduling and Dispatch}

The frequency-security constraints deduced in this paper can be used to solve the scheduling of a power grid while guaranteeing that enough inertia and FR will be available in the event of an outage. After the commitment solution is fixed, and therefore the inertia level is decided, the system operator could then send the order to the chosen providers for frequency response to adjust their frequency control: in the case of droop-controlled generators, the droop gains would be changed to guarantee the agreed FR deliverable by the specified time after a fault; in the case of converter-based devices (such as battery storage or smart loads), their fast power injections are similar to ramps controlled by activation signals, so their FR dynamics would closely match the ones considered in eq. (\ref{definitionStaticFR}) without further tuning. 

Before solving the frequency-constrained scheduling, the system operator could gather information from all devices willing to provide frequency response: these providers would have to submit the fastest delivery time that the device can comply with under any possible system condition (i.e. any level of inertia combined with contingency size), which for the case of generators would likely come from their ramp limitations. With this information, the system operator would solve the frequency-secured UC, and after that would inform FR providers of the chosen level of inertia and contingency size so that their frequency controls can be tuned appropriately.

This proposed approach is consistent with current practice in Great Britain (and to the best of our knowledge, in other power systems in the world), which simply establishes that a particular service such as Primary Frequency Response must be fully delivered 10 seconds following a contingency. PFR providers have to guarantee delivery of the agreed amount of FR by this time, no matter the size of the contingency. 

While FR is expected to be increasingly delivered by converter-based devices, which therefore would entail no conservativeness in the dispatch of FR using the constraints proposed here, if in a given system most of the providers of FR are droop-controlled generators, future work could be oriented towards reducing conservativeness in the FR solution. That is, studying ways in which full potential of FR delivery can be extracted by more precisely mapping the ramp parameters considered in the frequency constraints ($\textrm{T}_s$  and $\textrm{T}_{\textrm{del},s}$) to the droop dynamics. This mapping would need to be done before solving the frequency-secured UC, when the inertia and contingency size are not yet known, so that the FR provider can submit appropriate dynamic parameters to the system operator before the frequency-secured scheduling is solved.

Finally, a note on the need for communication channels: the methodology presented here is agnostic as to how the activation signal for FR is sent. This can be done through local frequency measurements for the droop control of a generator or for the control system of a power-electronics-based device such as BESS, or through a communication channel used by an aggregator to trigger different distributed devices. All these cases can be modelled in the proposed formulation by considering appropriate activation delays for each FR service. More detailed discussion on ICT requirements for frequency response services can be found in \cite{QitengDelays}.

\subsection{Applicability to Markets for Ancillary Services}

Most electricity systems are nowadays based on markets for energy, and it is not uncommon that some ancillary services are also procured through auctions or similar mechanisms. While it is out of the scope of the present paper to discuss market arrangements for frequency services, which could in practice take several shapes including the simultaneous clearing with pool energy markets, the proposed frequency-constrained formulation could also be applied to a market for frequency services, particularly so given that the proposed formulation is convex.

This convex formulation allows to implement a marginal-pricing scheme for frequency services, proposed and described in detail in a partner paper \cite{LuisPricing}. However, the present paper focuses on developing analytical constraints and understanding the value for the system from optimizing the different frequency response services that could be extracted from available providers.

Focusing on the particular case of GB, an island where the low-inertia problem is more acute than in continental grids, the system operator National Grid is currently designing a new suite of frequency response products. Although the new services are still being designed, they are envisioned to be modelled in the same way as proposed in this paper, i.e. as a power-increase ramp in combination with an activation delay \cite{NationalGridDC_FAQ}. Furthermore, an auction to competitively procure inertia was held at the beginning of 2020 \cite{NationalGridInertiaTender}. Therefore, the frequency-secured formulation presented here would allow the system operator to optimally clear the ancillary services auctions, instead of pre-defining a volume of response to be procured through each separate auction, as is current practice.

\section{Case Studies} \label{SectionCaseStudies}

\begin{table}[!t]
\renewcommand{\arraystretch}{1.2}
\caption{Characteristics of Thermal Plants in the GB 2030 System}
\label{TableThermal}
\centering
\begin{tabular}{l| l l l l}
    \multicolumn{1}{c|}{} & Nuclear & CCGT & OCGT\\
\hline
Number of Units & 4 & 100 & 30\\
Rated Power (MW) & 1800 & 500 & 100\\
Min Stable Generation (MW) & 1400 & 250 & 50\\
No-Load Cost (\pounds/h) & 0 & 4500 & 3000\\
Marginal Cost (\pounds/MWh) & 10 & 47 & 200\\
Startup Cost (\pounds) & N/A & 10000 & 0\\
Startup Time (h) & N/A & 4 & 0\\
Min Up Time (h) & N/A & 4 & 0\\
Min Down Time (h) & N/A & 1 & 0\\
Inertia Constant (s) & 5 & 4 & 4\\
Max FR deliverable (MW) & 0 & 50 & 20\\
\hline
\end{tabular}
\end{table}

In order to highlight the importance of co-optimizing the provision of distinct FR services, as well as to understand the value of different services under diverse system conditions, several case studies were carried out in a Stochastic Unit Commitment model \cite{AlexEfficient}, with the frequency-security constraints deduced in Section \ref{SectionFrequency} implemented. This SUC considers the uncertainty in wind forecast and optimally schedules energy and reserve leading to significant operating savings in low carbon systems, as demonstrated in \cite{AlexEfficient}; by adding post-fault frequency constraints, the solution of the SUC is also guaranteed to maintain frequency stability. 

The characteristics of the system under consideration are given in Table \ref{TableThermal}. A 10GWh pump-storage unit is included, with 2.6GW rating and 75\% round efficiency, corresponding to the Dinorwig plant in GB. Battery Energy Storage Systems (BESS) with 90\% efficiency and a 5h tank are also present, with a capacity of 200MW. 

Simulations spanning one year of operation were run, with frequency-security requirements of $\textrm{RoCoF}_{\textrm{max}}=0.5\textrm{Hz/s}$ and $\Delta f_{\textrm{max}}=0.8\textrm{Hz}$, while $\textrm{P}_{\textrm{L}}^{\textrm{max}}=1.8\textrm{GW}$. The quantiles for the SUC were set to 0.005, 0.1, 0.3, 0.5, 0.7, 0.9 and 0.995.  

In the frequency-secured scheduling model considered here there is no explicit cost for providing frequency services, i.e. there is no explicit term in the objective function corresponding to inertia or FR. The only costs in the optimization are fuel costs from thermal generators. Since frequency response is mainly provided by keeping some headroom in thermal generators, and system inertia is increased if a higher number of generators are committed, both these services imply running part-loaded generators, which increases generation costs as thermal units operate at a less efficient point below rated power and additional wind generation may be curtailed.  

The cost reduction from co-optimizing faster FR services shown in some of the case studies below is therefore driven by the lower amount of total inertia and FR needed to secure the nadir. In other words, less part-loaded plants are needed online if faster FR is optimized, therefore decreasing the operational cost of the system.

\subsection{Importance of Defining and Co-Optimizing New FR Services} \label{SectionNewFR}

\begin{figure} [!t]
    \centering
    \includegraphics[width=3.3in]{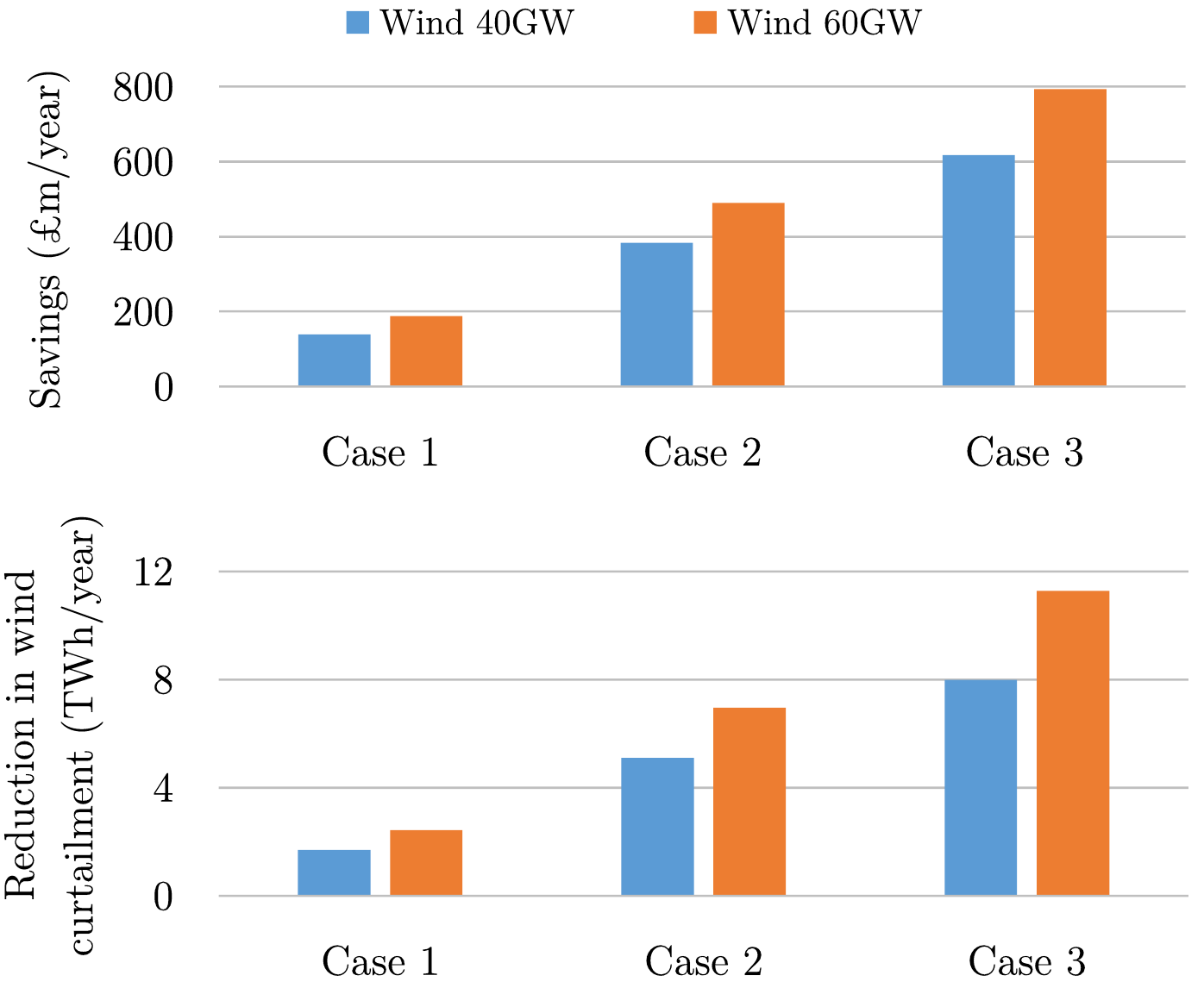}
    \caption{Benefits, in terms of economic savings and reduction in wind curtailment, from defining new FR services as compared to the state-of-the-art methods \cite{VincenzoEFR,LuisEFR}.}
	\label{Importance_fig1}
\end{figure}

Creating new FR services involves a trade-off between improved market efficiency and increased market complexity. A fundamental question that needs to be answered is ``How many and which new FR services should be defined?" While previous models \cite{VincenzoEFR,LuisEFR} only allow to co-optimize up to two distinct FR services, here the benefits of optimizing additional FR services by applying the proposed frequency-secured formulation are shown, as this formulation can efficiently co-optimize any number of FR services.
This section therefore focuses on quantifying the value in defining new FR services using the GB 2030 system, while the results give insight on the benefits that new services would bring to any different system. 

The frequency-secured optimization proposed in this paper is the first one allowing to co-optimize multiple FR services, with any combination of activation delays for each service. The state-of-the-art methods \cite{VincenzoEFR} and \cite{LuisEFR} allow to optimize up to two FR services, without giving the flexibility to consider different activation delays. Other works in the literature such as \cite{FeiStochastic}, \cite{OPFChavez} and \cite{LinearizedUC} consider a single FR service, and the relative benefits from incrementally considering two services have been studied in \cite{VincenzoEFR} and \cite{LuisEFR}. Therefore, the case studies in the present paper are focused on understanding the incremental benefits from co-optimizing three and four FR services, with respect to the two FR services in \cite{VincenzoEFR} and \cite{LuisEFR}.

In this section, some Combined Cycle Gas Turbines (CCGTs) are assumed to provide FR in less than 10s, corresponding to a ``fast PFR" service that some gas plants could provide \cite{NationalGridEFCC}. 
The generation mix in Table \ref{TableThermal} is considered, and from the total number of CCGTs, 70\% are assumed to provide PFR in 10s, 20\% have the capability of providing FR in 7s and the 10\% remaining have the capability of providing FR in 5s. Four different cases for FR services are defined: 
\begin{itemize}
 \item Base case, corresponding to the state-of-the-art \cite{VincenzoEFR,LuisEFR}: only EFR and PFR are defined by the system operator, as is current practice in GB. Therefore, all CCGTs are considered to provide PFR even if some of them can actually achieve faster FR dynamics. EFR is provided by the BESS.
 \item Case 1: a new FR service is defined, $\textrm{FR}_2$, delivered in 7s (i.e. $\textrm{T}_2=7\textrm{s}$), therefore the CCGTs with the capability of providing FR in 5s and 7s can provide this $\textrm{FR}_2$.
 \item Case 2: a new FR service $\textrm{FR}_2$ is defined  with $\textrm{T}_2=5\textrm{s}$, therefore the CCGTs with the capability of providing FR in 5s can provide this $\textrm{FR}_2$, while the CCGTs with the capability of providing FR in 7s will provide PFR.
 \item Case 3: two new FR services are defined, $\textrm{FR}_2$ with $\textrm{T}_2=5\textrm{s}$ and $\textrm{FR}_3$ with $\textrm{T}_3=7\textrm{s}$. This allows to fully extract the value of the different dynamics in FR provision available.
\end{itemize}
For each of these cases, two different wind-capacity levels are considered: 40GW and 60GW. 
The nuclear units are assumed to be fully loaded, therefore $P_{\textrm{L}}=\textrm{P}_{\textrm{L}}^{\textrm{max}}=1.8\textrm{GW}$ .

The results are presented in Fig. \ref{Importance_fig1}, showing the benefits of the three cases referred to the base case. It is interesting to note that Case 2 shows higher benefits than Case 1: by having only 10\% of the CCGTs providing FR in 5s in Case 2 (all other CCGTs are assumed to provide PFR in 10s), higher savings can be achieved than in Case 1, where 30\% of the CCGTs provide FR in 7s (20\% with the capability of providing FR in 7s plus 10\% with the capability of 5s). The results clearly demonstrate the complex task of defining new services, which will be a trade-off between the FR speed of delivery and the amount of provision. Finally, Case 3 shows the benefits from taking full advantage of fast FR from CCGTs, by defining two new FR services with delivery times of 5s and 7s. For all cases, the savings increase with wind penetration, as higher wind capacity implies lower inertia available from thermal generators, and therefore recognizing fast FR services becomes more valuable. 

Note that the benefits of co-optimizing fast FR services are not only in terms of cost, but also in reduced wind curtailment, as shown in Fig. \ref{Importance_fig1}: by defining the new FR services, a lower number of thermal generators need to stay online simply for providing inertia and FR, and therefore more wind power can be accommodated. The wind-curtailment reduction in Fig. \ref{Importance_fig1} is again referred to the base case, for which wind curtailment was of 31.86TWh/year for the 40GW-wind-capacity scenario (which means a 26.36\% of wind energy curtailed) and of 79.97TWh/year for the 60GW-wind-capacity (44.10\% wind energy curtailed).

\subsection{Computational Performance}
\begin{table}[!t]
\renewcommand{\arraystretch}{1.9}
\caption{Computation time for one year of system operation}
\label{table_computation}
\centering
\begin{tabular}{|>{\centering\arraybackslash}m{3.0cm}|>{\centering\arraybackslash}m{1.3cm}|>{\centering\arraybackslash}m{1.3cm}|>{\centering\arraybackslash}m{1.9cm}}
    \multicolumn{1}{c|}{} & Cases 1\&2 & Case 3 & 10 FR services \\ 
\hline
\multicolumn{1}{c|}{Single thread} & 2h 50min & 5h 10min & 89h 20min \\
\hline
\multicolumn{1}{c|}{Four threads} & 20min & 55min & 30h 20min \\
\hline

\end{tabular}
\end{table}

The simulation time for the cases presented in Section \ref{SectionNewFR} is included in Table \ref{table_computation}. Note that these simulations represent a whole year of operation of the GB 2030 system, for the computationally intensive Stochastic UC problem. The SUC optimizations were solved with FICO Xpress 8.0 in a 3.5GHz Intel Xeon CPU with twelve cores and 64GB of RAM, and the duality gap for the MISOCPs was set to 0.5\%. The multi-thread capabilities of the SUC framework used \cite{AlexEfficient} have been exploited to demonstrate the reduction in computational time that can be achieved through multi-core computing.

To further analyse the computational burden of the frequency-secured SUC problem, an additional case is considered in Table \ref{table_computation}, in which 10 FR services are defined. This case has been created using Case 3 as a base: the generators considered to provide FR by 5s in Case 3 are split into three new FR services with delivery times of 4s, 4.5s and 5s. The generators considered to provide FR by 7s are split into new services of 5s, 5.5s, 6s, 6.5s and 7s. Finally, the generators providing FR by 10s are split into services of 8.5s and 10s, making the total number of FR services defined of 10.

This case with 10 FR services has been added for illustrative purposes, as two main difficulties are foreseen for defining 10 different services in the operation of a real power system: 1) the fleet of FR providers must be significantly heterogeneous so as to have so distinct delivery times; and 2) the market complexity could significantly increase for both the system operator and market participants to guarantee delivery of FR within such specific times (note that some FR services are in this case only 0.5s apart from the consecutive service). Nevertheless, if any particular system does indeed have such heterogeneous dynamics for FR, and the non-synchronous capacity is high, this example demonstrates that quite significant extra savings could be achieved by defining even more FR services: savings from the 10 FR services case are £999m/year from the Base case, for 60GW wind penetration. This implies further £205m/year savings when compared to the £794m/year savings reported for Case 3 in Fig. \ref{Importance_fig1}. The drawback is a significant increase in computation complexity, which can however be reduced by taking advantage of multi-core computing.

The increase computation time when defining new FR services is related to the additional binary decision variables needed for the conditional statements in the nadir constraints (\ref{ConditionalSOCconstraints}). Note that the computational performance reported in Table \ref{table_computation} corresponds to a full year of operation of the system using the stochastic rolling planning approach, where each time-step in the SUC optimization solves a 24h look-ahead horizon problem. For the most complex case considered, the 10 FR services, each time-step was solved on average in 37s using a single thread, which would make it suitable for a real-time system dispatch.

Regarding the performance of previously proposed frequency-secured scheduling methods, reference \cite{VincenzoEFR} reports computational times ranging between 0.9h to 2.75h, for a case in which only two different FR services were considered. Note however that the comparison of computational performance of the methodology proposed in the present paper with that on \cite{VincenzoEFR} cannot be made on a one-to-one basis: computational time depends mostly on the characteristics of the optimization not related to FR, that is, if the UC solved is deterministic or stochastic, if integer variables are relaxed, and the granularity of the rolling planning in the UC.

\subsection{Impact of the Availability of Frequency Services on the Benefits from Defining New Services}

\begin{figure} [!t]
    \centering
    \includegraphics[width=3.3in]{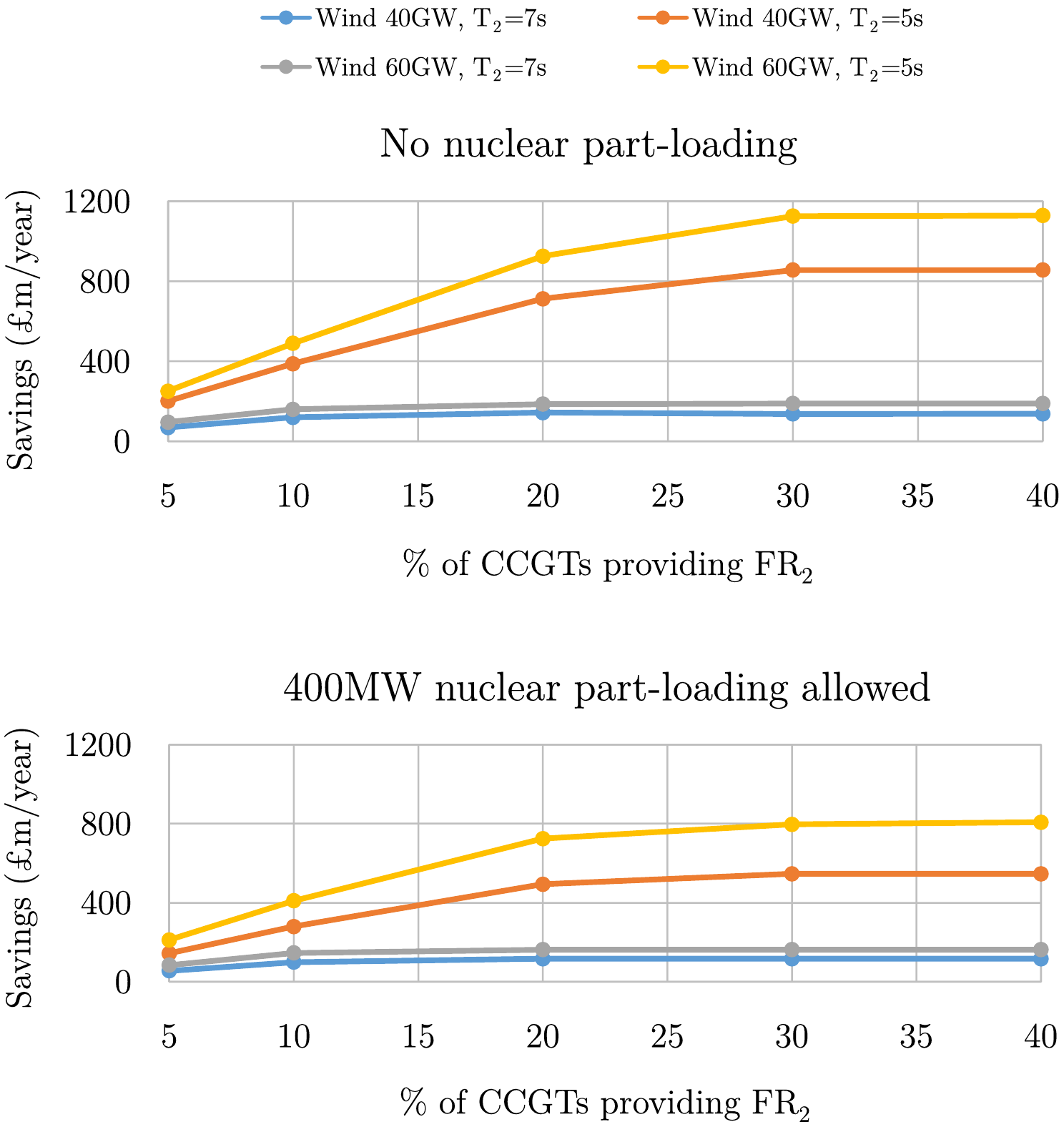}
    \caption{Impact of the mix of providers of frequency-services on the savings from defining a new FR service.}
	\label{FigImpactMix}
\end{figure}

The previous section has demonstrated that the value of new FR services will depend on the speed of delivery of these services, 
but another factor must be taken into account: the availability of frequency services, i.e. the number of units/providers that can deliver each FR service.
On a practical note, the available devices willing to provide FR would communicate it to the system operator ahead of real-time delivery of the service, as the current practice in frequency response markets.

In order to analyze such impact, a new FR service is considered to have been created, and the following variations are studied: 1) the new $\textrm{FR}_2$ service is provided by 5\%, 10\%, 20\%, 30\% and 40\% of the total CCGTs; and 2) same as the previous case, but nuclear units are allowed to part-load 400MW from their rated power (note that nuclear part-loading is co-optimized along with every other frequency service). Furthermore, wind penetrations of 40GW and 60GW are considered, as well as delivery times for $\textrm{FR}_2$ of $\textrm{T}_2=7\textrm{s}$ and $\textrm{T}_2=5\textrm{s}$. The ``base case" is the same as in Section \ref{SectionNewFR}, where all CCGTs provide PFR.

The results are presented in Fig. \ref{FigImpactMix}, which shows that the benefits from defining a new $\textrm{FR}_2$ service highly depend on the mix of providers: the savings can be limited if only 5\% of the CCGTs can provide $\textrm{FR}_2$, particularly if there is already 400MW of nuclear part-loading available and $\textrm{T}_2=7\textrm{s}$ (savings of \pounds55m/year for a 40GW-wind scenario). The first case titled ``No nuclear part-loading" shows increasing savings with respect to the base case when more units can deliver the new $\textrm{FR}_2$ service, but a saturation effect is also present: when the new $\textrm{FR}_2$ service defined has a $\textrm{T}_2=7s$, going beyond 20\% of the CCGTs providing this $\textrm{FR}_2$ does not further increase savings, while the saturation occurs at 30\% of the CCGTs if the new $\textrm{FR}_2$ service is delivered with $\textrm{T}_2=5s$. Moreover, the results for the second case titled ``400MW nuclear part-loading" show reduced savings for all scenarios of percentage of CCGTs providing $\textrm{FR}_2$, wind penetration level, and delivery time for $\textrm{FR}_2$: this is due to a reduced need for fast FR services when the output level of the largest units (the nuclear units for the GB system) is co-optimized along with FR. In other words, there would exist competition between the ``part-loading nuclear" service and fast FR.

In conclusion, these results demonstrate that the system operator should conduct a survey to check how many providers can potentially provide new fast FR services, before starting the process to change the market rules to define such FR services.

\subsection{Impact of Activation Delays in FR Provision} \label{SectionDynamicStatic}

\begin{figure} [!t]
    \centering
    \includegraphics[width=3.3in]{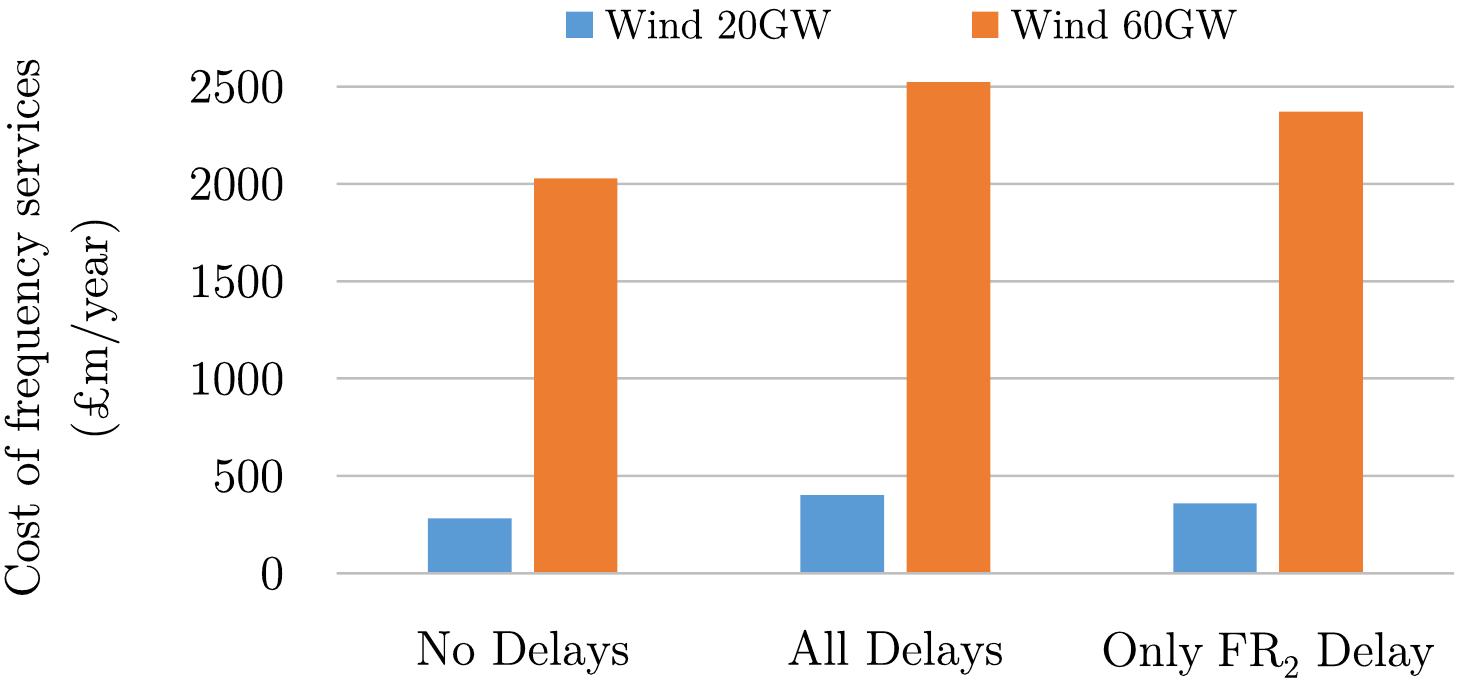}
    \caption{Impact of activation delays of FR services on the operational cost of the system.}
	\label{FigCaseStudyDelays}
\end{figure}

\begin{figure} [!t]
    \centering
    \includegraphics[width=3.3in]{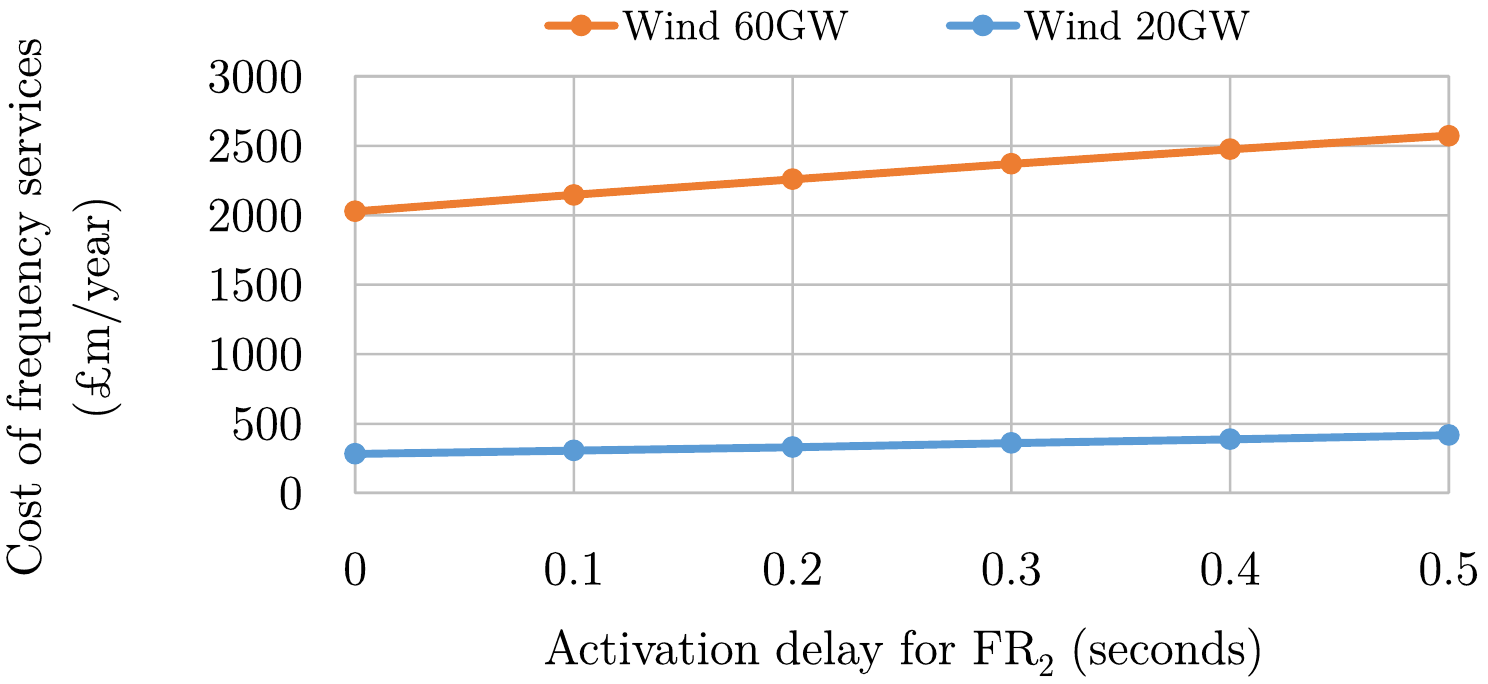}
    \caption{Sensitivity analysis for the impact of the time-delay in service $\textrm{FR}_2$ on the operational cost of the system.}
	\label{FigCaseStudyDelays2}
\end{figure}

This section analyzes the impact on system operating cost due to activation delays for the different FR services, i.e. the time after the generation outage when the FR services are activated and start ramping up. The base case here assumes 30\% of the CCGTs providing $\textrm{FR}_2$ with $\textrm{T}_2=5\textrm{s}$. Two different wind capacities are studied: 20GW and 60GW. Then, three cases of FR delays are considered:
1) none of the FR services have a delay, therefore they all react to any deviation from nominal frequency; 2) all three FR services have a 0.3s delay; and 3) EFR and PFR have no delay, but $\textrm{FR}_2$ has a 0.3s delay. While the impact of FR delays on the frequency nadir has been studied through dynamic simulations in \cite{QitengDelays}, no previous work has studied their impact on the system operating cost.

The results are presented in Fig. \ref{FigCaseStudyDelays}, which shows the annual cost of frequency services for each case, that is, the cost related to providing inertia and the different types of FR. This ``cost of frequency services" is calculated by taking the cost from the solution of the frequency-secured SUC minus the cost of an SUC with no frequency constraints. By comparing the second case ``All Delays" to the first one ``No Delays", it is demonstrated that delays in the delivery of FR services can significantly reduce their value to the system. However, a 0s delay would mean reacting to any deviation from nominal frequency, which would likely have associated wear-and-tear for the device providing FR. This is particularly true for battery storage providing EFR, since the lifetime of the device can be greatly impacted from the frequent charge-discharge switching that would be caused by reacting to any frequency deviation. Finally, the third case ``Only $\textrm{FR}_2$ delay" shows that by eliminating the delay in the provision of EFR and PFR, the system costs can be reduced by more than \pounds150m/year for a 60GW-wind scenario. The economic impact of time-delays for FR provision shown in Fig. \ref{FigCaseStudyDelays} is very significant for the 60GW-wind scenario, but this impact is limited for the 20GW-wind scenario.

Finally, Fig. \ref{FigCaseStudyDelays2} presents a sensitivity analysis for the activation delay in $\textrm{FR}_2$, for the case ``Only $\textrm{FR}_2$ delay", demonstrating that a reduction of just 0.1s in $\textrm{T}_{\textrm{del},2}$ can have a significant impact in annual system costs: reducing the delay from 0.3s to 0.2s would bring \pounds111m/year savings, for the 60GW-wind scenario.

\subsection{Role of Forecast for Inertia Contribution from Demand}

\begin{figure} [!t]
    \centering
    \includegraphics[width=3.3in]{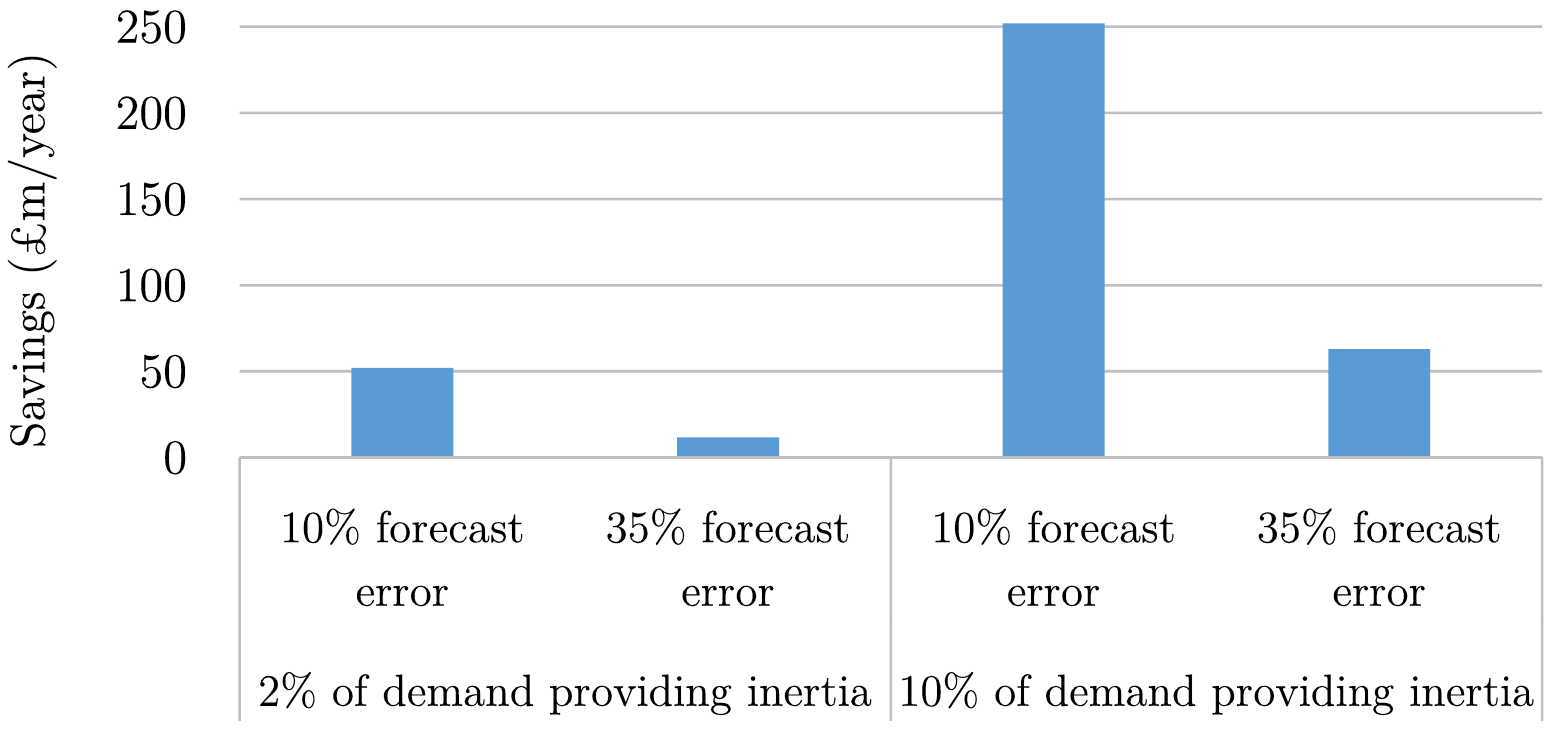}
     \vspace{1pt}
    \caption{Savings from considering the inertia contribution from demand, under different percentages of demand providing inertia and different forecast errors.}
	\label{FigInertiaDemand}
    \vspace{4pt}
\end{figure}

\begin{figure} [!t]
    \centering
    \includegraphics[width=3.3in]{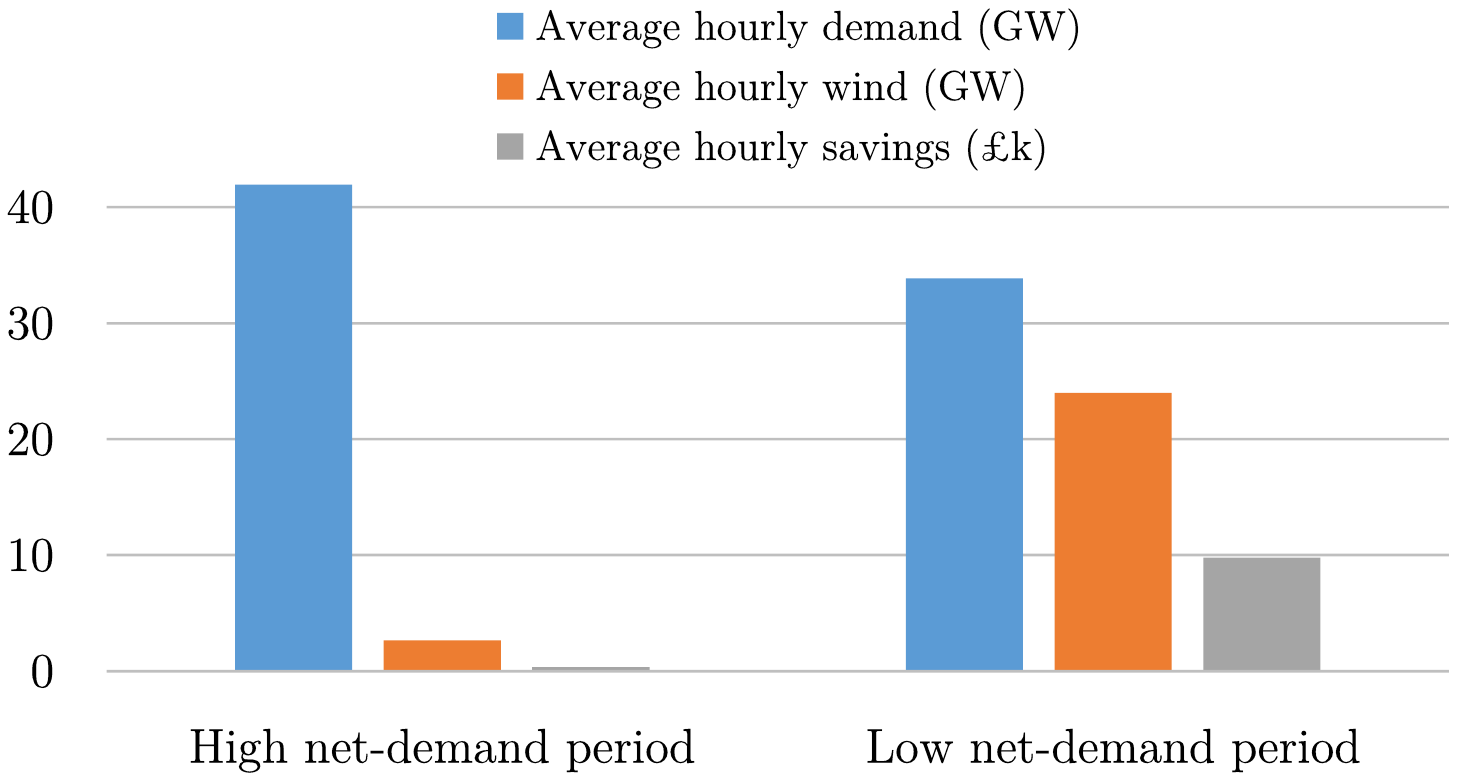}
    \vspace{1pt}
    \caption{Average hourly savings due to considering the inertia contribution from the demand side, for two cases corresponding each to a 50-hour period.}
	\label{FigInertiaDemand2}
\end{figure}

In this section the system operating cost is studied, under different forecasting scenarios for the inertia contribution from the demand side. The base case here considers 30\% of the CCGTs providing $\textrm{FR}_2$ with $\textrm{T}_2=7\textrm{s}$ and a 40GW wind capacity. The system operator is assumed to require a 99\% probability for the RoCoF and nadir constraints (\ref{reformulationRocofChance}) and (\ref{reformulationNadirChance}) to be fulfilled, i.e. $\alpha=\eta=0.99$. 

From this base case, in which no inertia from demand is taken into account, four further cases are considered: 2\% and 10\% of the total demand at each time-step in the SUC is assumed to provide inertia, with a forecast error of 10\% and 35\% considered for each case (i.e. $\sigma=0.1\textrm{H}_\mu$ or $\sigma=0.35\textrm{H}_\mu$). 
The inertia constant for demand is assumed to be 5s in all cases. The results in Fig. \ref{FigInertiaDemand} demonstrate that considering the inertia contribution from demand can bring non-negligible economic savings to the system, more so if a higher percentage of the demand is contributing to inertia. Furthermore, the results highlight the importance of achieving an accurate forecast for the inertia contribution from demand: by reducing the forecast error from 35\% to 10\%, \pounds190m/year additional savings can be obtained, for a case of 10\% of demand contributing to inertia.

Furthermore, Fig. \ref{FigInertiaDemand2} shows the detailed operation of the system by considering two cases, each corresponding to a 50-hour period, but showing significantly different levels of net-demand (demand minus wind): a high net-demand period in January, and a low net-demand period in July.
The results in Fig. \ref{FigInertiaDemand2}, which consider that 2\% of the demand provides inertia and the forecast error is of 10\%, demonstrate how valuable inertia from demand becomes during low net-demand periods: the savings in the second period (average hourly savings of \pounds9.8k) are significantly higher than in the first period (average hourly savings of \pounds0.36k), although there is less inertia available from the demand side in the former. 
The results clearly demonstrate that market arrangements need to be in place not only to incentivize the demand side to provide inertia, but also to incentivize these inertia-providing loads to consume during periods when low net-demand is expected.

\section{Conclusion and Future Work} \label{SectionConclusion}
This paper proposes frequency-security constraints that allow to consider any finite number of FR providers, as well as any combination of activation delays for the different FR services. Uncertainty in system inertia from the demand side is modelled using chance constraints, for which a convex reformulation is provided. The resulting MISOCP formulation allows to efficiently schedule FR considering the different dynamics of providers, co-optimizing FR along with inertia and a reduced largest loss. Case studies of frequency-secured Stochastic Unit Commitment have demonstrated the importance of new FR services for achieving cost-efficiency in power systems with significant renewable penetration.

Regarding future work, the reliability of FR provision should be considered, since it has been assumed here that all FR providers can deliver the agreed amount of FR. Inertia from the demand side has been shown to provide important savings, and therefore new methods to accurately forecast this inertia should be developed. Finally, the different regional frequencies in certain buses of the power grid should be modelled, as this work has considered the uniform-frequency model driven by the Centre of Inertia.

%
%
%
%



\ifCLASSOPTIONcaptionsoff
  \newpage
\fi



%



\IEEEtriggeratref{19}

\bibliographystyle{IEEEtran} 
\bibliography{Luis_PhD}

\begin{thebibliography}{10}
\providecommand{\url}[1]{#1}
\csname url@samestyle\endcsname
\providecommand{\newblock}{\relax}
\providecommand{\bibinfo}[2]{#2}
\providecommand{\BIBentrySTDinterwordspacing}{\spaceskip=0pt\relax}
\providecommand{\BIBentryALTinterwordstretchfactor}{4}
\providecommand{\BIBentryALTinterwordspacing}{\spaceskip=\fontdimen2\font plus
\BIBentryALTinterwordstretchfactor\fontdimen3\font minus
  \fontdimen4\font\relax}
\providecommand{\BIBforeignlanguage}[2]{{%
\expandafter\ifx\csname l@#1\endcsname\relax
\typeout{** WARNING: IEEEtran.bst: No hyphenation pattern has been}%
\typeout{** loaded for the language `#1'. Using the pattern for}%
\typeout{** the default language instead.}%
\else
\language=\csname l@#1\endcsname
\fi
#2}}
\providecommand{\BIBdecl}{\relax}
\BIBdecl

\bibitem{FeiISGT2017}
F.~Teng, M.~Aunedi, G.~Strbac, V.~Trovato, and A.~Dallagi, ``Provision of
  ancillary services in future low-carbon {UK} electricity system,'' in
  \emph{2017 IEEE PES Innovative Smart Grid Technologies Conference Europe
  (ISGT-Europe)}, Conference Proceedings.

\bibitem{ElaI}
E.~Ela, V.~Gevorgian, A.~Tuohy, B.~Kirby, M.~Milligan, and M.~O'Malley,
  ``Market designs for the primary frequency response ancillary service. {P}art
  {I}: Motivation and design,'' \emph{IEEE Transactions on Power Systems},
  vol.~29, no.~1, pp. 421--431, 2014.

\bibitem{ERCOT_EFR}
W.~Li, P.~Du, and N.~Lu, ``Design of a new primary frequency control market for
  hosting frequency response reserve offers from both generators and loads,''
  \emph{IEEE Transactions on Smart Grid}, vol.~9, no.~5, pp. 4883--4892, 2018.

\bibitem{EDF_cuts}
C.~Cardozo, W.~{van Ackooij}, and L.~Capely, ``Cutting plane approaches for
  frequency constrained economic dispatch problems,'' \emph{Electric Power
  Systems Research}, vol. 156, pp. 54--63, 2018.

\bibitem{OPFChavez}
H.~Ch\'{a}vez, R.~Baldick, and S.~Sharma, ``Governor rate-constrained {OPF} for
  primary frequency control adequacy,'' \emph{IEEE Transactions on Power
  Systems}, vol.~29, no.~3, pp. 1473--1480, 2014.

\bibitem{IowaThesis}
G.~Zhang, ``New ancillary service market design to improve {MW}-frequency
  performance: reserve adequacy and resource flexibility,'' Ph{D} thesis,
  {I}owa {S}tate {U}niversity, {U}nited {S}tates, 2015.

\bibitem{FeiStochastic}
F.~Teng, V.~Trovato, and G.~Strbac, ``Stochastic scheduling with
  inertia-dependent fast frequency response requirements,'' \emph{IEEE
  Transactions on Power Systems}, vol.~31, no.~2, pp. 1557--1566, 2016.

\bibitem{LinearizedUC}
H.~Ahmadi and H.~Ghasemi, ``Security-constrained unit commitment with
  linearized system frequency limit constraints,'' \emph{IEEE Transactions on
  Power Systems}, vol.~29, no.~4, pp. 1536--1545, 2014.

\bibitem{PricingElaZhang}
G.~Zhang, E.~Ela, and Q.~Wang, ``Market scheduling and pricing for primary and
  secondary frequency reserve,'' \emph{IEEE Transactions on Power Systems},
  vol.~34, no.~4, pp. 2914--2924, 2019.

\bibitem{VincenzoEFR}
V.~Trovato, A.~Bialecki, and A.~Dallagi, ``Unit commitment with
  inertia-dependent and multispeed allocation of frequency response services,''
  \emph{IEEE Transactions on Power Systems}, vol.~34, no.~2, pp. 1537--1548,
  2019.

\bibitem{LuisEFR}
L.~Badesa, F.~Teng, and G.~Strbac, ``Simultaneous scheduling of multiple
  frequency services in stochastic unit commitment,'' \emph{IEEE Transactions
  on Power Systems}, vol.~34, no.~5, pp. 3858--3868, 2019.

\bibitem{UCFaroe}
L.~E. Sokoler, P.~Vinter, R.~Bærentsen, K.~Edlund, and J.~B. J{\o}rgensen,
  ``Contingency-constrained unit commitment in meshed isolated power systems,''
  \emph{IEEE Transactions on Power Systems}, vol.~31, no.~5, pp. 3516--3526,
  2016.

\bibitem{NationalGridPLossInertia}
\BIBentryALTinterwordspacing
``System needs and product strategy,'' National Grid, Report, 2017. [Online].
  Available:
  \url{https://www.nationalgrid.com/sites/default/files/documents/8589940795-System%20Needs%20and%20Product%20Strategy%20-%20Final.pdf}
\BIBentrySTDinterwordspacing

\bibitem{QitengDelays}
Q.~Hong, M.~Nedd, S.~Norris, I.~Abdulhadi, M.~Karimi, V.~Terzija, B.~Marshall,
  K.~Bell, and C.~Booth, ``Fast frequency response for effective frequency
  control in power systems with low inertia,'' \emph{The Journal of
  Engineering}, vol. 2019, no.~16, pp. 1696--1702, 2019.

\bibitem{OMalleyDeload}
R.~Doherty, G.~Lalor, and M.~O'Malley, ``Frequency control in competitive
  electricity market dispatch,'' \emph{IEEE Transactions on Power Systems},
  vol.~20, no.~3, pp. 1588--1596, 2005.

\bibitem{LuisPESGM2018}
L.~Badesa, F.~Teng, and G.~Strbac, ``Optimal scheduling of frequency services
  considering a variable largest-power-infeed-loss,'' in \emph{2018 IEEE Power
  {\&} Energy Society General Meeting (PESGM)}, Conference Proceedings.

\bibitem{KundurBook}
P.~Kundur, \emph{Power System Stability and Control}, 1st~ed.\hskip 1em plus
  0.5em minus 0.4em\relax McGraw-Hill Education, 1994.

\bibitem{NGrocofDefinition}
\BIBentryALTinterwordspacing
``System operability framework,'' National Grid, Report, 2016. [Online].
  Available:
  \url{https://www.nationalgrid.com/sites/default/files/documents/8589937803-SOF%202016%20-%20Full%20Interactive%20Document.pdf}
\BIBentrySTDinterwordspacing

\bibitem{ConejoOptimizationBook}
R.~Sioshansi and A.~Conejo, \emph{Optimization in Engineering}.\hskip 1em plus
  0.5em minus 0.4em\relax Springer International Publishing, 2017.

\bibitem{BoydConvex}
S.~Boyd and L.~Vandenberghe, \emph{Convex Optimization}.\hskip 1em plus 0.5em
  minus 0.4em\relax Cambridge University Press, 2004.

\bibitem{SaphiroChance}
A.~Nemirovski and A.~Shapiro, ``Convex approximations of chance constrained
  programs,'' \emph{SIAM Journal on Optimization}, vol.~17, no.~4, pp.
  969--996, 2007.

\bibitem{RandomVariableBook}
A.~Papoulis and S.~U. Pillai, \emph{Probability, Random Variables and
  Stochastic Processes}, 4th~ed.\hskip 1em plus 0.5em minus 0.4em\relax
  McGraw-Hill Education, 2002.

\bibitem{LuisPricing}
L.~Badesa, F.~Teng, and G.~Strbac, ``Pricing inertia and frequency response
  with diverse dynamics in a {Mixed-Integer Second-Order Cone Programming}
  formulation,'' \emph{Applied Energy}, vol. 260, p. 114334, 2020.

\bibitem{NationalGridDC_FAQ}
\BIBentryALTinterwordspacing
``Dynamic {C}ontainment {FAQ}s,'' National Grid, Report, 2020. [Online].
  Available: \url{https://www.nationalgrideso.com/document/164391/download}
\BIBentrySTDinterwordspacing

\bibitem{NationalGridInertiaTender}
\BIBentryALTinterwordspacing
``Stability {P}athfinder for inertia, phase 1 tender results,'' National Grid,
  Report, 2020. [Online]. Available:
  \url{https://www.nationalgrideso.com/document/162091/download}
\BIBentrySTDinterwordspacing

\bibitem{AlexEfficient}
A.~Sturt and G.~Strbac, ``Efficient stochastic scheduling for simulation of
  wind-integrated power systems,'' \emph{IEEE Transactions on Power Systems},
  vol.~27, no.~1, pp. 323--334, 2012.

\bibitem{NationalGridEFCC}
\BIBentryALTinterwordspacing
``Enhanced frequency control capability project,'' National Grid, Report, 2017.
  [Online]. Available:
  \url{https://www.nationalgrideso.com/innovation/projects/enhanced-frequency-control-capability-efcc}
\BIBentrySTDinterwordspacing

\end{thebibliography}




\vspace{21.5mm}

\vskip -0.2\baselineskip plus -1fil

\begin{IEEEbiographynophoto}{Luis Badesa}
(S'14) received the B.S. in Industrial Engineering degree 
from the University of Zaragoza, Spain, in 2014, and the 
M.S. in Electrical Engineering degree from the University
of Maine, United States, in 2016. Currently he is pursuing
a Ph.D. in Electrical Engineering at Imperial College London, U.K. His research interests lie in modelling
and optimization for the operation of low-carbon power grids.
\end{IEEEbiographynophoto}

 

\vskip -0.5\baselineskip plus -1fil

\begin{IEEEbiographynophoto}{Fei Teng}
(M'15) received the BEng in Electrical Engineering
from Beihang University, China, in 2009, and the 
MSc and PhD degrees in Electrical Engineering from 
Imperial College London, U.K., in 2010 and 2015.
Currently, he is a Lecturer in the Department of Electrical 
and Electronic Engineering, Imperial College London, U.K.
His current research areas are the operation of low-inertia power systems, 
data analytics and cyber-physical power systems.
\end{IEEEbiographynophoto}

\vskip -0.5\baselineskip plus -1fil

\begin{IEEEbiographynophoto}{Goran Strbac}
(M'95) is Professor of Electrical Energy
Systems at Imperial College London, U.K. His current
research is focused on the optimisation of operation and investment of low-carbon energy systems, energy infrastructure reliability and future energy markets. 
\end{IEEEbiographynophoto}




\end{document}